\documentclass[11pt]{article}
\usepackage{graphicx}
\usepackage{amssymb}
\usepackage{epstopdf}
\newtheorem{theorem}{Theorem}
\newtheorem{cor}[theorem]{Corollary}
\newtheorem{lem}{Lemma}
\newtheorem{definition}{Definition}

\title{Deconstructing  Functions on Quadratic Surfaces into Multipoles}
\author{Gabriel Katz}

\begin{document}
%
%

\maketitle

\begin{abstract} Any homogeneous polynomial $P(x, y, z)$ of degree $d$, being restricted to a unit sphere $S^2$, admits essentially a unique representation of the form  $\lambda_0 + \sum_{k = 1}^d \lambda_k [\prod_{j = 1}^k L_{kj}]$, where $L_{kj}$'s are linear forms in $x, y$ and $z$ and $\lambda_k$'s are real numbers. The coefficients of these linear forms, viewed as 3D vectors, are called \emph{multipole} vectors of $P$. In this paper we consider similar multipole representations of polynomial and analytic functions on other quadratic surfaces $Q(x, y, z) = c$, real and complex. Over the complex numbers, the above representation is not unique, although the ambiguity is essentially finite. We investigate the combinatorics that depicts this ambiguity. We link these results with some classical theorems of harmonic analysis, theorems that describe decompositions of functions into sums of spherical harmonics. We extend these classical theorems (which rely on our understanding of the Laplace operator $\Delta_{S^2}$) to more general differential operators $\Delta_Q$ that are constructed with the help of the quadratic form $Q(x, y, z)$.  Then we introduce modular spaces  of multipoles. We study their intricate  geometry and topology using methods of algebraic geometry and singularity theory. The multipole spaces are ramified over vector or projective spaces, and the compliments to the ramification sets give rise to a rich family of $K(\pi, 1)$-spaces, where $\pi$ runs over a variety of modified braid groups.   
\end{abstract} 

\section{Introduction: the cosmic background motivation}

This paper aims to generalize and extend some results in  [KW] and [W] about deconstruction of cosmic microwave background radiation (CMBR) into multipole vectors and explain these results to a mathematically inclined audience. Our exposition, to a degree, is self-sufficient.

Cosmologists who study the CMBR routinely decompose it into 
a sum of spherical harmonics---the eigenfunctions of the Laplace operator on the 
$2$-sphere.  The summand that corresponds to the eigenvalue $-d(d + 1)$ of the Laplace operator  is called a $d$-th \emph{spherical harmonic}. This decomposition helps  to analyze the correlations between the magnitude and geometry of various  harmonics. The ultimate goal of these investigations is to understand the geometry of the visible universe and the physics of the Big Bang.  

The data obtained by Wilkinson Microwave Anisotropy Probe (WMAP) reveal puzzling correlations between the low-$d$ portion of the harmonic decomposition of the cosmic microwave background radiation [B], [CHS], [TOH], [EHGL]. In the case of quadrupole-octopole ($d = 2$ and $3$) alignment, simulations based on our fast algorithm [KW], [W] show the alignment to be unusual at the 98.7 percent level. No explanation is yet known for these strange results.

{\it Multipole vectors} provide a convenient means to study the isotropic properties of CMBR. In particular, they help to quantify the planarity of a given multipole, as well as to compare the
alignment of two different multipoles [CHS].  Jeff Weeks and I,
 coming from a background in pure mathematics, were unable
to decipher  the formalism and terminology of [CHS] and
chose instead to re-create the multipole vector concept from
scratch (see [KW]).  Later on, as we read the papers  [L] and [D], we realized that we managed 
to reinvent the wheel twice: our key algebraic lemma  turns out to be a classical 
theorem of Sylvester [S], [S1], while the application of that lemma to harmonic polynomials 
on sphere implicitly was known to Maxwell [M].  

The real-valued spherical harmonics of order $d$ are precisely the homogeneous harmonic polynomials of degree $d$ in the variables $x$, $y$ and $z$\footnote{For example, the spherical harmonic $Y_{2}^{0}$ is the polynomial $x^2 + y^2 - 2 z^2$, up to normalization.}. Thus, we sought to understand the multipole vectors of Copi, Huterer and Starkman [CHS] from a polynomial point of view. As in [CH], this calls for using some elementary tools from the field of algebraic geometry.  
Translated to the language of polynomials, [CHS]'s motivating goal was to factor every homogeneous
harmonic polynomial $P$ of degree $d$ on the sphere into a product of linear factors
\begin{eqnarray}
P(x,y,z) =  (a_1 x + b_1 y + c_1 z)  (a_2 x + b_2 y + c_2 z) 
\quad \dots \quad (a_d x + b_d y + c_d z).
\end{eqnarray}
Such a factorization is  impossible in general, as [CHS] implicitly acknowledge by their introduction of suitable error terms.  The correct statement  is given by the Sylvester Theorem [S], [S1]:\footnote{rediscovered  in [WK] by the  authors.}
\begin{theorem} 
Every real homogeneous polynomial $P$ of degree $d$ in $x$, $y$ and $z$ may be written as
\begin{eqnarray}
P(x,y,z) &=&  (a_1 x + b_1 y + c_1 z)  (a_2 x + b_2 y + c_2 z)       
\quad ... \quad (a_d x + b_d y + c_d z) + \nonumber\\
&+&  (x^2 + y^2 + z^2)\cdot R(x, y, z),
\end{eqnarray} 
where the remainder term $R$ is a homogeneous polynomial
of degree $d - 2$.  The decomposition is unique up to
reordering and rescaling of the linear factors in the product.
\end{theorem}
For a \emph{harmonic} homogeneous $P$, a somewhat similar representation was discovered by James Clerk Maxwell in his 1873 \emph{Treatise on Electricity  and Magnetism} [M].  Its relation to algebraic geometry in general, and to the Sylvester theorem, in particular, is  well-explained  by Courant and Hilbert in [CH] VII, Sec. 5. We will describe the Maxwell representation later in the paper. 

Sylvester's Theorem 1 (which is just an application of  B\'{e}zout's theorem) has a pleasing corollary:
\begin{theorem} { $[KW]$}
When restricted to a unit sphere $S^2$, every real polynomial $P$ in $x, y$, and $z$ of degree $d$ can be written in the form 
\begin{eqnarray}
P(x, y, z) = &
\lambda_0 + \lambda_1 (a_{11} x + b_{11} y + c_{11} z) + \lambda_2 (a_{21} x + b_{21} y + c_{21} z)(a_{22} x + b_{22} y + c_{22} z) + \quad \dots \nonumber \\
 & \dots \quad + \lambda_d  (a_{d1} x + b_{d1} y + c_{d1} z)(a_{d2} x + b_{d2} y + c_{d2} z) \; \dots \;
(a_{dd} x + b_{dd} y + c_{dd} z)
\end{eqnarray}
The decomposition is unique up to reordering and rescaling of the linear factors in each of the products.
In particular, one can pick each  $(a_{jk}, b_{jk},  c_{jk})$  to be a unit vector and, for  $j > 0$, each $\lambda_j \geq 0$.
\end{theorem}
It is well-known that spherical harmonics form a basis in the space $L_2 (S^2)$ of $L_2$-integrable  functions on a unit sphere (see [CH]). As a result, an infinite Maxwell-type decomposition is available for any function $f \in L_2 (S^2)$ (cf.  [L], formula $(26)$).
\begin{theorem}  Any real function $f \in L_2(S^2)$  has a representation in the form
\begin{eqnarray}
f(x, y, z) = 
\sum_{d = 1}^\infty \Big \{ \lambda_{d, 0} +  \sum_{k = 1}^d \lambda_{d, k}  \prod_{j = 1}^k (a_{dj} x + b_{dj} y + c_{dj} z) \Big \}
\end{eqnarray}
where  the series converges in the $L_2$-norm, and the mutually orthogonal polynomials in the figure brackets are among the $d$-th spherical harmonics.  This representation is unique,  up to reordering and rescaling of the linear factors in each of the products.

In the space of real continuous  functions $f$ on $S^2$ there is a dense (in the $sup$-norm) and $O(3)$-invariant subset $\mathcal F(S^2)$ comprised of functions that admit a representation in the form of a series
\begin{eqnarray}
 \lambda_0 + \sum_{d = 1}^\infty   \lambda_d  \prod_{k = 1}^d(a_{k} x + b_{k} y + c_{k} z) 
 \end{eqnarray}
 that uniformly converges  on the sphere. Each functions from $\mathcal F(S^2)$ admits a canonical  extension into the interior of the unit ball $D^3$ where it produces a real analytic function.  
\end{theorem}
This theorem, claiming a form of  stability of decomposition $(3)$ under the polynomial approximations, is a useful tool for analyzing patterns in the spherical sky [W].  

At the first glance, the sphere seems to occupy a special place in these results. However, what is really important, is the \emph{quadratic nature} of the surface.  This paper is concerned with  similar decompositions of functions on other quadratic algebraic skies (real and complex). For instance, it is natural to wonder:  Can one deconstruct in a similar fashion  polynomial or real analytic functions on a hyperboloid?  The short answer is ``Yes, one can", but the uniqueness of the representation is lost.

The focal point of this paper is to describe rich and interesting algebro-geometrical structures and  topology of  modular spaces of  multipoles. In the process, we bring under the same roof a variety of mathematical techniques and constructions that belong to the fields of algebraic geometry, singularity theory, algebraic topology, harmonic analysis, and  the polynomial approximation theory. None of our techniques is very advanced, but their  natural appearance within the context of studying quadratic surfaces is pleasing and somewhat surprising...

The reader can follow  two distinct treads in the core of  the paper. The first unifies results that are linear in nature and are based on the classical theory of harmonic analysis on quadratic surfaces. These results are valid for quadratic hypersurfaces \emph{in any dimension}. The second tread winds through the results that require methods of algebraic geometry. These results are uniquely  \emph{three-dimensional}, that is, applicable only to quadratic \emph{surfaces}.


\section{Multipoles and polynomials on complex quadratic surfaces}

We denote by $\mathbb C[x, y, z]$ the ring of complex polynomials in the variables $x, y,$ and $z$. Let 
$Q(x, y, z)$ be an irreducible quadratic form over the complex numbers $\mathbb C$. From a strictly  algebraic perspective, one can interpret  this section  as describing the ways in which the quotient ring 
$\mathbb C[x, y, z] /  \langle Q - \lambda \rangle$, 
$\lambda \in \mathbb C$, fails to be a Unique Factorization Domain. However, the flavor of our approach to this problem is more geometrical and combinatorial.

The reader will be well-adviced to ignore for a while the asterisks in our notations: they are there to distinguish between vector spaces and their duals. Let ${\mathcal{S} } = \{ Q(x, y, z) = 1\}$ be a complex algebraic surface in $\mathbb C^3_\ast$.  
At the same time, $\{ Q(x,y, z) = 0\}$ gives rise to a complex projective curve  $\mathcal{Q} \subset \mathbb CP^2_\ast$, where $\mathbb CP^2_\ast$ stands for the projectivization of the $3$-space with coordinates $x, y$, and $z$.
As in Theorem 2, we aim to decompose any polynomial $P(x, y, z)$  restricted to $\mathcal S$ as an ``economic"  sum of products of homogeneous linear polynomials.

First, consider the case of an homogeneous $P$ and the corresponding complex projective quadratic curve $\mathcal{P}$ in $\mathbb{C}P^2_\ast$ it generates.  Denote by $Z(P, Q)$ the intersection $\mathcal{P} \cap \mathcal{Q}$. We assume that  $\mathcal{P}$ and  $\mathcal{Q}$
do not share a common component, so that $Z(P, Q)$ is a finite set.  When  $Z(P, Q)$ is a complete intersection, it consists of exactly $2d$ points, where  $d = deg(P)$. 
\begin{definition} Let $Q(x, y, z)$ be an homogeneous irreducible quadratic polynomial, and $P(x, y, z)$ an homogeneous polynomial of degree $d$. Assume that $Z(P, Q)$ is a complete intersection.  A  \emph{parcelling} of the set $Z(P, Q)$ is comprized of a number of subsets $\{Z_\nu \subset Z(P, Q)\}_\nu$ so that: 
\begin{itemize}
\item
$Z(P, Q) = \coprod_{\nu} Z_\nu$ 
\item
 $Z_\nu \cap Z_{\nu'} = \emptyset$, provided $\nu \neq \nu'$
\item
the cardinality of each subset $Z_\nu$ is equal to 2.
\end{itemize}
\end{definition}
When $Q$ is irreducible and $P$ is not divisible by $Q$, the intersection set $Z(P, Q)$ is 
equipped  with a  function $\mu$ which assigns to each point $p \in Z(P, Q)$ its multiplicity $\mu(p)$. 

We denote by $\mathbb Z_+$ the set of non-negative integers and by   $\mathbb N$ the set of positive integers.
\begin{definition} Let  $Z$ be a finite set equipped with a multiplicity function 
$\mu: Z \rightarrow \mathbb N$ whose $l_1$-norm $\|\mu\|_1$ is 
$2d$.  A \emph{generalized} \emph{parcelling} of $(Z, \mu)$ is a collection of functions 
$\mu_\nu : Z \rightarrow \{0, 1, 2\}$, such that 
\begin{itemize}
\item $\sum_\nu \mu_\nu = \mu$
\item $\|\mu_\nu\|_1 = 2$
\end{itemize}
\end{definition}
Of course, each parcelling of $Z(P, Q)$ is also a generalized parcelling, where 
the roles of the functions $\mu_\nu$ are played by the characteristic functions of the parcels $Z_\nu$. 

In the case of $Q = x^2 + y^2 + z^2$, the following lemma is due to Sylvester [S], [S1].
\begin{lem} Let $Q(x, y, z)$ be an irreducible homogeneous quadratic polynomial  and let $P(x, y, z)$ be
a homogeneous polynomial of degree $d$ which is not divisible by $Q$.   Consider  a generalized  parcelling $\mu = \sum_\nu \mu_\nu$ 
of the multiplicity function $\mu : Z(P, Q) \rightarrow \mathbb N$.
 
Then the homogeneous polynomial $P(x, y, z)$  admits a  representation of the form 
\begin{eqnarray}
Q(x, y, z)\cdot R(x, y, z) + \prod_\nu L_\nu (x, y, z) ,
\end{eqnarray}
where  $L_\nu $ denotes a linear homogeneous polynomial  that  vanishes at each point 
$p \in Z(P, Q)\subset \mathcal Q$ with the multiplicity $\mu_\nu(p)$. 
\end{lem}
{ \bf Proof \quad}  
There exist a linear polynomial $L_\nu$  and a corresponding line  
$\mathcal L_\nu \subset  \mathbb{C}P^2_\ast$  such that the multiplicity of $\mathcal L_\nu \cap \mathcal Q$ at each point $p \in Z(P, L_\nu)$ equals $\mu_\nu(p)$.  When $\mu_\nu(p) = 2$, the line must be tangent to $\mathcal Q$ at $p$. When the support of $\mu_\nu$---the parcel 
$Z_\nu \subset Z(P, Q)$)---is comprised of two points, the line $\mathcal L_\nu$ is chosen to contain $Z_\nu$. 

Let us compare the restrictions of $P$ and $\prod_\nu L_\nu$ to the curve $\mathcal Q$.
Both polynomials are of the same degree $d$. Moreover, 
the curve $\mathcal L := \{\prod_\nu L_\nu = 0\}$ intersects with $\mathcal Q$ at  the points of 
$Z(P, Q)$, where it realizes the  multiplicities prescribed by the function  $\mu$.
It follows that the restrictions of $P$ and $\prod_\nu L_\nu$ to $\mathcal Q$ are proportional, that is,   for an appropriate choice of scalar $\lambda$,\, 
$P|_{\mathcal Q} = \lambda \cdot \prod_\nu {L_\nu}|_{\mathcal Q }$. Just take $\lambda = P(q)/\prod_\nu L_\nu(q)$ where $q \in \mathcal Q \setminus Z(P, Q)$. With this choice, the curves 
$\{P - \lambda \prod_\nu {L_\nu} = 0\}$ and $\mathcal Q$ intersect so that the total multiplicity of the intersection is at least $2d + 1$. 
By the B\'{e}zout theorem, this is possible only when 
$ Z(P - \lambda \prod_\nu {L_\nu})  \supset \mathcal Q$. 
Employing the irreducibility of $Q$, we get that 
$P - \lambda \cdot  \prod_\nu L_\nu$ must be divisible by $Q$.
This completes the proof. $\Box$
\begin{cor}  Any effective divisor $D$  on the complex quadratic surface $\mathcal S = \{Q(x, y, z) = 0\}$ that is  the zero set of a homogeneous polynomial  $P(x, y, z)$ can be represented as a sum of lines $D_\nu$. $\Box$
\end{cor}
In general, the representation (6) of $P$ depends on a generalized parcelling, and thus, is not unique. 
\smallskip

We denote by $V(d)$ the complex vector space of homogeneous  polynomials in $x, y$, and $z$ of degree $d$. Its dimension is $(d^2 + 3d + 2)/2$. Consider a vector subspace 
$V_Q(d) \subset V(d)$ comprised of polynomials divisible by $Q$. There is a canonical imbedding 
$\beta_Q : V(d - 2) \rightarrow V(d)$ whose image is $V_Q(d)$. It is produced by multiplying  polynomials of degree $d - 2$ by $Q$. Thus, any space $V(d)$ is equipped with a natural filtration 
$\mathcal F_Q(d) = \{V_{Q^k}(d)\}_k$ by subspaces of polynomials divisible by various powers $\{Q^k\}$ of $Q$.
We  equip each $V(d)$ with an Hermitian inner product. For various $d$'s, we insist  that these inner products are  synchronized by the requirement: all the imbeddings $\beta_Q$  must be isometries.
Denote by $V_Q^\perp(d)$ the subspace of $V(d)$ orthogonal to $V_Q(d)$.  It can be identified with  the quotient  space $V(d)/ V_Q(d)$---the $d$-graded part of the algebra of regular functions on the complex cone $\{Q(x, y, z) = 0\}$. The  dimension of $V_Q^\perp(d)$ is equal to $(2d + 1)$. 

For a real form $Q$, similar spaces based on real homogeneous polynomials of a degree $d$ make sense. We denote them $V(d; \mathbb R), V_Q(d; \mathbb R)$, and $V_Q^\perp(d; \mathbb R)$.

 When $Q = x^2 + y^2 + z^2$, there is a classical interpretation for the quotient  
 $V(d; \mathbb R)/ V_Q(d; \mathbb R)$. It can be identified with the set of \emph{harmonic} homogeneous polynomials of degree $d$ (cf.  [CH] and [M]). In other words, the kernel  
$Har(d; \mathbb R)$ of the Laplace operator 
$\Delta: V(d; \mathbb R) \rightarrow V(d - 2; \mathbb R)$  is complementary to 
$V_Q(d; \mathbb R)$  in $V(d; \mathbb R)$.
This leads to a decomposition similar to the one in $(6)$ (cf.  [Sh, Theorem 22.2], [L], and [W]) :
\begin{eqnarray} 
V(d; \mathbb R) \approx Har(d; \mathbb R) \oplus V_Q(d; \mathbb R)
\end{eqnarray}
The direct summands in $(7)$ are orthogonal with respect to the inner product in $V(d; \mathbb R)$ defined by the formula $\langle f, g \rangle = \int_{S^2} f \cdot g \; dm$. The measure  $dm$ on the sphere $S^2$  is the standard one.

Moreover, according to Maxwell, any  polynomial $P \in Har(d; \mathbb R)$ admits a beautiful representation of the form
\begin{eqnarray}
P(x, y, z) = r^{2d + 1} \nabla_{\mathbf v_1}\nabla_{\mathbf v_2} ... \nabla_{\mathbf v_d} \Big(\frac{1}{r}\Big),
\end{eqnarray}
where $r = (x^2 + y^2 + z^2)^{1/2}$, $\{\mathbf v_j \in \mathbb R^3\}$ are some  
vectors,\footnote{Later,  we  call them "the leading multipole vectors" of $P$.} and $\nabla_{\mathbf v_j}$ stands for the directional derivative operator. 

Physics behind Maxwell's representation is quite transparent: $1/r$ is a potential of a single electrical charge, $\nabla_{\mathbf v_1}(1/r)$ is a potential of a ``virtual" (that is, very small) dipole formed by two opposite charges, $\nabla_{\mathbf v_2}\nabla_{\mathbf v_1}(1/r)$ is a potential of a ``virtual" quadropole---a close pair of $\mathbf v_1$-oriented dipoles merging along the direction of 
$\mathbf v_2$---, $\nabla_{\mathbf v_3}\nabla_{\mathbf v_2}\nabla_{\mathbf v_1}(1/r)$ is a  potential of a ``virtual" octopole, and so on...

Note that, the Laplace operator $\Delta = \partial_x^2 + \partial_y^2 + \partial_z^2$ can act on complex polynomials as well. Moreover, if $\Delta(f) = 0$ for some  analytic function $f$,  then its real and imaginary parts are also harmonic:
$\Delta(\mathrm{Re} f) = 0$ and $\Delta(\mathrm{Im} f) = 0$. 

For any invertible complex $(3\times 3)$-matrix $A = (a_{j k})$,  consider a  symmetric matrix 
$B = (b_{j k}) = A\cdot A^T$ and the corresponding quadratic  form 
$Q(v) = \langle vB, \; v\rangle$, where $v = (x,  y, z)$,  $\langle\sim , \sim \rangle$ stands for the  inner product $x x' + y y' + z z'$ in $\mathbb C^3$, and the upper script $T$ denotes the matrix transposition.
Employing $Q$, we can form  a second order  differential operator 
$\Delta_Q = \sum_{1 \leq i, j \leq 3} b^{j k} \partial_j \partial_k$ acting on holomorphic 
functions $f$ in the complex variables $\{x_1 = x,\,  x_2 = y, \, x_3 = z\}$. Here $\{b^{ij}\}$ denote elements of the inverse matrix $B^{-1}$. Formally, $\Delta_Q = [\vec{\partial}]\cdot B^{-1} \cdot   [\vec{\partial}]^T$ where $[\vec{\partial}] := (\partial_1, \partial_2, \partial_3) := (\partial_x, \partial_y, \partial_z)$.

Given a non-degenerated quadratic form $Q$, there is a change of complex coordinates 
$(x', y', z') =  (x, y, z) \cdot A$ that reduces it to the canonical form $Q' = x'^2 + y'^2 + z'^2$. Consider the complex  Laplace operator 
$\Delta' = \partial_{x'}^2 +\partial_{y'}^2 + \partial_{z'}^2$ in the new coordinates $(x', y', z')$.
Then, for any holomorphic function $f(x', y',  z')$, we have 
$$[\Delta'  f]( (x, y, z)\cdot A) = \Delta_Q [f ((x, y, z)\cdot A)]$$ 
Thus, there is a 1-to-1 correspondence between harmonic  homogeneous polynomials 
$f(x', y', z')$ of degree $d$ and  degree $d$ homogeneous polynomial solutions 
$g(x, y, z)$ of the equation $$\Delta_Q(g(x, y, z)) = 0.$$
Let us examine the intersection $Ker(\Delta') \cap V_{Q'}(d)$. For any homogeneous polynomial $T$ of degree $d - 2$, we get $\Delta'(Q'\cdot T) = Q'\cdot \Delta'(T)  + \Delta'(Q')\cdot T + 2 \langle\nabla T, \nabla Q'\rangle = Q'\cdot \Delta'(T)  + (4d - 6)T$. Therefore, if $Q'\cdot T \in Ker(\Delta')$, then $T$ must be divisible by $Q'$ (note that $\Delta'(T) = 0$  implies $T = 0$). Put $T = Q'\cdot T_1$ and $\kappa(d) = 4d + 2$.  The equation 
$0 = \Delta'(Q'\cdot T) = Q'\cdot \Delta'(Q'\cdot T_1)  +  \kappa(d - 2) \cdot Q'\cdot T_1$ is equivalent to the equation 
$0 =  \Delta'(Q'\cdot T_1)  + \kappa(d - 2)\cdot T_1 = Q'\cdot \Delta'(T_1)  + [\kappa(d - 2) + \kappa(d - 4)]\cdot T_1$. Again, it follows that $T_1$ must be divisible by $Q'$. Continuing inductively this kind of reasoning, we see that $Ker(\Delta') \cap V_{Q'}(d) = \{0\}$. On the other hand, one can verify that 
$dim [Ker(\Delta')]  + dim [V_{Q'}(d)] = dim [V(d)]$. Thus, $V(d) = Ker(\Delta') \oplus V_{Q'}(d)$. 

A polynomial $P(x', y', z')$ is divisible by $x'^2 + y'^2 + z'^2$, if and only if, $P((x, y, z)\cdot A)$ is divisible by $Q(x, y, z)$. Therefore, $V(d) = Ker(\Delta_Q) \oplus V_{Q}(d)$ as well.

Let $O_Q(3; \mathbb C)$ denote a subgroup of the general linear group $GL(3; \mathbb C)$ that preserves the quadratic form $Q$. A matrix $U \in O_Q(3; \mathbb C)$ if and only if $UBU^T = B$. The natural $O_Q(3; \mathbb C)$-action on the space of $x, y$, and $z$-variables  induces an action  on the polynomial space $V(d)$. Evidently, $V_Q(d)$ is invariant under this action. 
On the other hand, for any 
polynomial $P(x, y, z)$,\, $$\Delta_Q [ P((x, y, z)\cdot U)] =  [\Delta_{\tilde Q} P]((x, y, z)\cdot U)),$$ where the operator  $\Delta_{\tilde Q} := [\vec{\partial}]\cdot U^T B^{-1} U \cdot   [\vec{\partial}]^T$.  Since 
$UBU^T = B$, we get $(U^{-1})^T B^{-1} U^{-1} = B^{-1}$.  By a simple algebraic trick, it follows that $U^T B^{-1} U = B^{-1}$. As a result, both the quadratic form $Q$ and the kernel $Ker(\Delta_Q)$ are invariant under the $O_Q(3; \mathbb C)$-action.

Consider Maxwell's representation (8) of a real homogeneous harmonic polynomial $P(x', y', z')$ of degree $d$. It gives rise to a map $\Xi$ that takes the sets of vectors 
$\mathbf v_1, \mathbf v_2, \dots,  \mathbf v_d$ to elements of $Har(d, \mathbb R)$.  The map  $\Xi$ is evidently linear in each of the $\mathbf v_j$'s, and thus it is a real polynomial map with a vector space of real dimension $2d + 1$ for its target. By [CH], [L] and [W], $\Xi$ is onto the vector space 
$Har(d, \mathbb R)$. Therefore, the complexification $\Xi^\mathbb C$ of $\Xi$ must be also onto the vector space 
$Har(d, \mathbb C) \approx Har(d, \mathbb R) \oplus \mathbf i \; Har(d, \mathbb R)$  (indeed, the image of $\Xi^\mathbb C$ in $Har(d, \mathbb C)$ is a complex algebraic set containing a totally real vector subspace $Har(d, \mathbb R)$ of a maximal dimension). In other words, formula (8) must be valid for any complex homogeneous harmonic polynomial $P(x', y', z')$ of degree $d$ and  appropriate complex vectors $\mathbf u_1, \mathbf u_2, \dots,  \mathbf u_d \in \mathbb C^3$. In fact, each vector $\mathbf u_j = \mathbf v_j \cdot A^{-1}$.
Formula (8)  describes any  homogeneous complex polynomial solution of the equation $\Delta_Q(f) = 0$ in the $(x', y', z')$ coordinates.  Translating them back to the $(x, y, z)$-coordinates with the help of the identity 
$[\nabla_{\mathbf v'} f'](\mathbf x\cdot A) = \nabla_{\mathbf v' \cdot A^{-1}}[ f'(\mathbf x \cdot A)]$  leads to a formula $(10)$ below. 
Thus, we have established the following proposition: 
\begin{lem}
The  decomposition
\begin{eqnarray}  
V(d) = Ker(\Delta_Q) \oplus V_Q(d)
\end{eqnarray}
holds for any non-degenerated complex quadratic form $Q$. It is invariant under the natural
$O_Q(3; \mathbb C)$-action on $V(d)$. Moreover, 
for any $P \in Ker(\Delta_Q)$ of degree $d$, there exist  vectors $\{\mathbf u_k  \in \mathbb C^3\}_{1 \leq k \leq d}$ so that the generalized Maxwell formula
\begin{eqnarray}  
P(x, y, z) = Q(x, y, z)^{d+ \frac{ 1}{2}} \cdot \nabla_{\mathbf u_1}\nabla_{\mathbf u_2} \dots 
\nabla_{\mathbf u_d}  \Big(Q(x, y, z)^{-\frac{1}{2}}\Big)
 \end{eqnarray}
 is valid. The real part of the LHS of $(10)$ generates all real homogeneous polynomial solutions $P$ of the equation $\Delta_Q(P) = 0$.

In particular, with $Q = x^2 + y^2 - z^2$, any  degree $d$ homogeneous complex polynomial solution $P$ of the wave  equation
 $(\partial^2_x + \partial^2_y)P = \partial^2_z P$    is given by the formula $(10)$. $\Box$
 \end{lem}
Formula $(10)$ might pose a slight challenge: after all, square roots are multivalued analytic functions. However, the $\pm$-ambiguities in picking a single-valued branch cancel each other. We shall see that, even up to reordering and rescaling of the multipole vectors $\{\mathbf u_k\}$, the complex Maxwell representation $(10)$ of a given ``$Q$-harmonic" $P$ is not unique. 
 \smallskip
 
 {\bf Remark.} In fact, by a similar argument, the decomposition (9) is available for homogeneous polynomials in \emph{any number of variables} and for any non-degenerated  quadratic form $Q$ (see [Sh], Theorem 22.2, for the proof). At the same time, the representation (10) seems to be a 3-dimensional phenomenon. It looks like that not any $Q$-harmonic polynomial in $n > 3$ variables can be expressed in terms of $d$ directional derivatives sequentially applied to the $Q$-harmonic potential $Q^{1 - n/2}$. This seems to be related to the failure of the Sylvester-type formula (6) for $n > 3$. When I raised this issue with Michael Shubin, he proposed a nice conjecture  in the flavor of Maxwell's representation (although, not a direct generalization of it). In the Maxwell  representation, one employs a differential operator which is a \emph{monomial} in $d$ directional derivatives, while in the conjecture below one invokes  a differential operator which is a \emph{polynomial} in $d$ directional derivatives.
 
 {\bf Shubin's Conjecture} ({\it An $n$-dimensional variation on the theme of  Maxwell's representation}). \hfill\break
 Let $n > 2$.  For any real  homogeneous and harmonic polynomial $P$ of degree $d$ in $n$ variables $\{x_j\}$,  there exists a unique real homogeneous and harmonic polynomial $P^\clubsuit$ also in $n$ variables such that 
 $$P(x_1, \dots ,  x_n)  = r^{2d + n -2} \Big[P^\clubsuit(\partial_{x_1}, \dots ,  \partial_{x_n})\; 
 r^{2 - n}\Big],$$
 where $r = (\sum_{j = 1}^n x_j^2)^{1/2}$, and  the differential operator 
 $P^\clubsuit(\partial_{x_1}, \dots ,  \partial_{x_n})$ being applied to the potential function $r^{2 - n}$. Moreover, the polynomials $P$ and  $P^\clubsuit$ are proportional with the coefficient of proportionality depending only on $d$ and $n$. $\Box$
 
Given a vector space $V$, we denote by $V^\circ$ the space $V \setminus \{0\}$. 
We consider  the $3$-space $V(1) \approx {\mathbb C^3}$ of linear forms $L = ax + by + cz$ (which we identify with the space of their coefficients $\{(a, b, c)\}$)  and the $3$-space $V_\ast(1) \approx \mathbb C^3_\ast$ with the coordinates $x, y, z$. This  calls for a distinction in notations for their projectivizations: $\mathbb CP^2$ and $\mathbb CP^2_\ast$. Hence, each point $(a, b, c) \in V(1)^\circ$ determines a point $l = [a : b : c]$ in $\mathbb CP^2$ and a complex line $\mathcal L$ in $\mathbb CP^2_\ast$--- the zero set of the form $L$.
\smallskip

Now, let us return to the decomposition $(6)$.
Each polynomial $L_\nu$ from $(6)$ can be viewed as a vector in $V(1) = V_Q^\perp(1)$.
Therefore,  a given homogeneous polynomial $P$ of degree $d$ which is not divisible by $Q$, together with the appropriate generalized parcelling, produces a collection of non-zero vectors 
$\{\mathbf{w}_\nu  \in V(1)\}_{\nu}$\footnote{These vectors  are not necasarilly  distinct.}.  
We shall call this unordered set of vectors $\{\mathbf{w}_\nu\}$ \emph{the leading multipole vectors of} $P$ with respect to the corresponding generalized parcelling  of $\mu: Z(P, Q) \rightarrow \mathbb Z_+$. 
\begin{lem} For a given generalized parcelling $\mu = \sum_\nu \mu_\nu$ of $\mu: Z(P, Q) \rightarrow \mathbb Z_+$, the subordinate  leading multipole vectors
$\{\mathbf{w}_\nu \in  V(1)\}_{\nu }$  are unique up to reordering and rescaling of the $L_\nu$'s (equivalently, the polynomial $R$ in $(6)$ is unique).
\end{lem}
{\bf  Proof}\quad We  use the same notations as in the proof of Lemma 1. By an argument as in that proof, any linear polynomial $L_\nu'$  which defines a line $\mathcal L_\nu'$  with the property 
$\mathcal  L_\nu'  \cap \mathcal Q = \mathcal  L_\nu \cap \mathcal Q$ (the points of intersection have multiplicities prescribed by $\mu_\nu$)  is of the form $\lambda L_\nu$.  Here $\lambda \in \mathbb C^\ast$ and  $L_\nu$ is a preferred  polynomial corresponding to a vector $\mathbf{w}_\nu \in V(1)$. 
Indeed, as in the proof of Lemma $1$,  for an appropriate $\lambda$, by the B\'{e}zout Theorem, $\lambda L_\nu - L_\nu'$ must be divisible by $Q$.  $\Box$

As a result, the multipole is well defined by the parcelling of $\mu$ modulo some rescaling.  Clearly, one can always replace each $L_\nu$ in $(6)$ with $\lambda_\nu L_\nu$ as long as 
$\prod_\nu \lambda_\nu = 1$. Next, we analyze the exact meaning of the ``reordering and rescaling" ambiguity.

Let $H_q$ be an abelian subgroup of $(\mathbb C^\ast)^q$ formed by vectors
$\{\lambda_\nu\}_{1 \leq \nu \leq q}$ subject to the restriction $\prod_\nu \lambda_\nu = 1$. Its rank is 
$q - 1$. Let $S_q$ stand for the symmetric group in $q$ symbols $\{1, 2, \dots , q\}$. We denote by $\Sigma_q$ an extension $1 \rightarrow H_q \rightarrow \Sigma_q \rightarrow S_q \rightarrow 1$ of $S_q$ by $H_q$. This group is generated by the obvious actions of  $S_q$ and $H_q$ on $(\mathbb C)^q$.

We introduce an orbit-space
\begin{eqnarray} 
\mathcal M(k) :=  [(V(1)^\circ]^k/ \Sigma_k
\end{eqnarray} 
whose points encode the products of linear forms as  in the decomposition $(6)$.
Here the group $\Sigma_k$ acts  on the products of  spaces by 
permuting them and by rescaling their vectors. Its subgroup $H_k$ acts freely. 

Because of  the uniqueness of the prime factorization in the polynomial ring $\mathbb C[x, y, z]$, 
the space  $\mathcal M(k)$ provides us with a $1$-to-$1$ parametrization of the space of homogeneous degree $k$ polynomials that are products of linear forms.
As an abstract space, (11)  can be expressed as
\begin{eqnarray}
 \mathcal M(k) =  
\Big \{[\mathbb C^{3 \, \circ}]^k \Big\}\Big/\Sigma_k
 \end{eqnarray} 
\begin{definition} The elements of  orbit-space $\mathcal M(k)$  will be  called $k$-\emph{poles}, or simply, \emph{multipoles}. 
\end{definition}
We introduce  groups $\Gamma_k$  in a manner similar to the introduction of the groups $\Sigma_k$.   The group $\Gamma_k$ is an extension of the permutation group $S_k$ by the group 
$(\mathbb C^\ast)^k \supset H_k$. Thus, $\Gamma_k /\Sigma_k  \approx \mathbb C^\ast$.

Therefore, the space $\mathcal M(k)$ in $(11), (12)$  fibers over the projective variety
\begin{eqnarray}
&\mathcal B(k) : = \prod_{\nu = 1}^k V(1)^\circ / \Gamma_k \nonumber\\  
& = Sym^k\mathbb CP[V(1)] \approx
 Sym^k (\mathbb CP^2)
\end{eqnarray} 
with the fiber $\mathbb C^\ast$.  Here $Sym^t(X) := X^t/S_n$ denotes the $t$-th symmetric power  of a space $X$.
The natural map $\mathcal M(k) \rightarrow \mathcal B(k)$ is a principle $\mathbb C^\ast$-fibration which gives rise to a line bundle 
$\eta(k) : = \{\mathcal M(k)\times_{\mathbb C^\ast}\mathbb C \rightarrow \mathcal B(k)\}$. 
By shrinking its zero section to a point, we form a quotient space 
 \begin{eqnarray}
 \overline {\mathcal M}(k) := [\mathcal M(k)\times_{\mathbb C^\ast}\mathbb C] \, /  \mathcal B(k).
 \end{eqnarray} 
It differs from $\mathcal M(k)$ by a single new point $\bf 0$---a point which will represent the zero multipole. Evidently, $\overline {\mathcal M}(k)$ is a contractible space.
 
Given a collection of vector spaces $\{V_\alpha\}$, their \emph{wedge} product $\wedge_\alpha V_\alpha$ (not to be mixed with the exterior product!) is the quotient of the Cartesian product $\times_\alpha V_\alpha$ by the subspace comprised of sequences $\{v_\alpha \in  V_\alpha\}$ such that at least one vector from the sequence is zero. 
Thus, topologically,
\begin{eqnarray} 
\overline{\mathcal M}(k) :=  [\wedge^{k}(V(1)]/ \Sigma_k.
\end{eqnarray}
Lemma 3 leads to the following proposition. \smallskip
\begin{cor} Consider an homogeneous  polynomial $P(x, y, z)$ of degree $d$ which is not divisible by an irreducible homogeneous quadratic polynomial $Q(x, y, z)$.  
Then any generalized parcelling $\mu = \sum_\nu \mu_\nu$ of the multiplicity function 
$\mu: Z(P, Q) \rightarrow \mathbb Z_+$  uniquely determines a multipole in the  space 
$\mathcal M(d)$ introduced in $(11)$ or $(12)$.  $\Box$
 \end{cor}
The space of multipoles has singularities. Its singular set $sing(\mathcal M(k))$ arises from the sets  points in $[\mathbb C^{3 \circ}]^k$ fixed by various non-trivial subgroups of $\Sigma_k$. These subgroups all are the conjugates (in $\Sigma_k$) of certain non-trivial subgroups of $S_k$ (recall that 
$H_k \subset \Sigma_k$ acts freely). The partially ordered set of the orbit-types  give rise to a natural stratification of the multipole space $\mathcal M(k)$. The space $sing(\mathcal M(k))$ is of complex codimension two in $\mathcal M(k)$. Its top strata corresponds to transpositions from $S_k$. Thus, a generic point from $sing(\mathcal M(k))$ has a normal slice in $\mathcal M(k)$ which is diffeomorphic to a  cone over the real projective space 
$S^3/\mathbb Z_2 = \mathbb RP^3$. The larger stabilizers of vectors from the space $\mathbb C^k$ of  the obvious $S_k$-representation  correspond to smaller strata of more complex geometry. 
Evidently, $sing(\mathcal M(k))$ is invariant under the diagonal action of 
$\mathbb C^\ast \approx \Gamma_k / \Sigma_k$.
These observations are summarized in the lemma below.
\begin{lem} The singular set $sing(\mathcal M(k)) \subset \mathcal M(k)$ is of codimension two. It is invariant under the $\mathbb C^\ast$-action.  Therefore, it is a principle $\mathbb C^\ast$-fibration over the singular set $sing(Sym^k(\mathbb CP^2)) \subset Sym^k(\mathbb CP^2)$.   A generic point of $sing(\mathcal M(k))$ has  $\mathbb RP^3$ for its normal link. 
Points of $sing(\mathcal M(k))$ correspond to  completely factorable polynomials 
$L = \prod_\nu L_\nu$ with at least two proportional linear factors (in other words, to  weighted collections of lines in $\mathbb CP^2_\ast$ that contain  at least one line of multiplicity  greater than one).  $\Box$
\end{lem}
Let ${\mathcal Fact}(d) \subset V(d)$ denote the variety of homogeneous polynomials of degree $d$ in $x, y$, and $z$ that are products of linear forms.

Given a multipole $w \in \mathcal M(d)$ one can construct the corresponding completely factorable polynomial $L(w) = \prod_\nu L_\nu \in V(d)$.  Note that, due to the uniqueness of the prime factorization in the polynomial ring and in view of our definition of multipoles, the correspondence  
$w \Rightarrow L(w)$ gives rise to a $1$-to-$1$ map  
$\Theta:  \mathcal M(d) \stackrel{\approx}{\rightarrow}   {\mathcal Fact}(d)$. 
Consider  a  ``Vi\`{e}te-type" algebraic map 
\begin{eqnarray} 
\Phi_Q: \mathcal M(d) \stackrel{\Theta}{\rightarrow}  {\mathcal Fact}(d)  \stackrel{\Pi_Q}{\rightarrow}  {\mathcal Fact}_Q(d),
 \end{eqnarray}
where $\Pi_Q$ is induced by restricting polynomials in $x, y$, and $z$ to the surface $\{Q(x, y, z) = 0\}$. The symbol  ${\mathcal Fact}_Q(d) \subset V(d) / V_Q(d) \approx V_Q^\perp(d)$ denotes the variety of homogeneous polynomial functions on the surface $\{Q = 0\}$ that also decompose into products of linear forms. Due to formula $(6)$ in Lemma 1, any non-zero homogeneous polynomial on the surface $\{Q = 0\}$ admits a linear factorization. Hence,  $\Phi_Q$ is \emph{onto} and ${\mathcal Fact}_Q(d)$ can be identified with the space $[V(d) / V_Q(d)]^\circ  \approx V_Q^{\perp}(d)^\circ$. 

The map $\Phi_Q$,  extends to an algebraic  map 
\begin{eqnarray}
\tilde \Phi_Q :  E \eta (d) \longrightarrow V_Q^{\perp}(d)
\end{eqnarray} 
defined on the space $E \eta (d)$ of the line bundle $\eta (d)$.  It sends the zero section 
$\mathcal B(d)$ of  $\eta (d)$ to the zero vector $\mathbf {0} \in  V_Q^{\perp}(d)$ and each 
fiber of $\eta (d)$ isomorphically to a line in $V_Q^{\perp}(d)$ passing through the origin. In fact,
$\hat \Phi_Q |_{E \eta (d)\setminus \mathcal B (d)} = \Phi_Q$. Evidently, $\hat \Phi_Q$ gives rise to a 
continuous map 
\begin{eqnarray}
\overline \Phi_Q : \overline \mathcal M(d) \longrightarrow  V_Q^{\perp}(d)
\end{eqnarray}
It  turns out that $\overline \Phi_Q$ has \emph{finite} fibers. We need some combinatorial constructions which will help us to prove this claim and to describe the cardinality of the $\overline \Phi_Q$-fibers. 

With every natural $d$   we associate an integer $\kappa(d)$ that counts the number of distinct parcellings  in a finite set of cardinality $2d$. 
Any parcelling is obtained by breaking a set $Z$ of cardinality $2d$ into  disjoint subsets of cardinality 2. Thus, $\kappa(d) = (2d - 1)!! : = (2d - 1)(2d - 3)(2d - 5) \quad... \quad 3\cdot 1$ is the number of possible handshakes among a company of $2d$ friends.  There is an alternative way to compute this number: consider the standard action of the permutation group $S_{2d}$ on the set of $2d$ elements. The action induces a transitive action on the set of all parcellings. Under this action, the  subgroup $S_{2d}^\#$ that preserves the parcelling  $\{\{1, 2\}, \{3, 4\},\; ...\;, \{2d -1, 2d\}\}$ is an extension 
\begin{eqnarray}
1 \rightarrow (S_2)^d \rightarrow S_{2d}^\#  \rightarrow S_d  \rightarrow 1
\end{eqnarray}
of the permutation group $S_d$ that acts naturally on the pairs by the group $(S_2)^d$ that exchanges  the elements in each pair. As a result, we get an  identity 
$$\kappa(d) = (2d - 1)!! = (2d)! / (2^d \cdot d!)$$  
As we deform a polynomial $P$ into a polynomial $P_1$, two or more points in $Z(P, Q)$ can merge into a single point of $Z(P_1, Q)$. Its multiplicity is equal to the sum of multiplicities of the points forming the merging group. Through this process, any generalized parcelling $p$ of the original 
$\mu: Z(P, Q) \rightarrow \mathbb N$ gives rise to a new and unique generalized parcelling $p_1$ of 
$\mu_1: Z(P_1, Q) \rightarrow \mathbb N$. Thus, we can define a partial order in the set of all generalized parcellings of  effective divisors  on $\mathcal Q$ of degree $2d$ by setting 
$\mu \succ \mu_1$ and $p \succ p_1$ (see Figures 1 and 3).

In the same spirit, let $\kappa(\mu)$ stand for the number of distinct generalized parcellings of a function $\mu: Z \rightarrow \mathbb N$ (recall that $\|\mu\|_1 = 2d$) on a finite set $Z$. Unless $\mu$ is identically $1$ and $|Z| = 2d$,  $\kappa(\mu) < \kappa(d)$.  When two intersection points 
merge, the number of generalized parcellings drops:  there are distinct original parcellings 
that  become indistinguishable after the merge (see the left diagram in Figure 1). For example, when two simple points in a complete intersection merge, the number of generalized parcellings changes from $\kappa(d)$ to 
$\kappa(d - 2) + [ \kappa(d) - \kappa(d - 2)]/2 =  [ \kappa(d) + \kappa(d - 2)]/2$, that is, it drops by 
$[ \kappa(d) - \kappa(d - 2)]/2$. In general, $\mu \succ \mu_1$ implies $\kappa(\mu) > \kappa(\mu_1)$. 

For a generic $L \in {\mathcal Fact}(d)$  the set $Z(L, Q)$, as well as all the parcels $\{Z_\nu\}$ defined by the linear factors $L_\nu$, are  complete intersections. For such an $L$, the number of multipoles in $\Phi_Q^{-1}(\Pi_Q(L))$ 
is the number $\kappa(d)$ of distinct parcellings in the set $Z(L, Q)$. For any 
$L \in {\mathcal Fact}(d)$, the cardinality of $\Phi_Q^{-1}(\Pi_Q(L))$ (equivalently, of 
$\Pi_Q^{-1}(\Pi_Q(L))$) equals to  the the number of generalized parcellings of the multiplicity function
$\mu: Z(L, Q) \rightarrow \mathbb N$. Indeed, assume that two completely factorable polynomials 
$L, L'$ coincide when restricted to the surface $\{Q = 0\}$. Then $L - L'$ must be divisible by 
$Q$. Therefore, $Z(L, Q) = Z(L', Q)$, moreover, the two multiplicities  of each point in the intersection (defined by the curves $\mathcal L$ and $\mathcal L'$) must be equal as well. Hence, $\mathcal L$ and $\mathcal L'$ define two parcellings of the same multiplicity function $\mu$ on the intersection set. 

Let $X, Y$ be topological spaces and  $f: X \rightarrow Y$ a continuous map with finite fibers. 
For a while,  the \emph{ramification set} $\mathcal D(f)$ of $f$ is understood as the set $\{y_0 \in Y\}$ such that, for  any open neighborhood  $U$ of $y_0$, the cardinality of the fibers $\{f^{-1}(y)\}_{y \in U}$ is not constant. 

\begin{figure}[ht]
\centerline{\includegraphics[height=1.5in,width=2in]{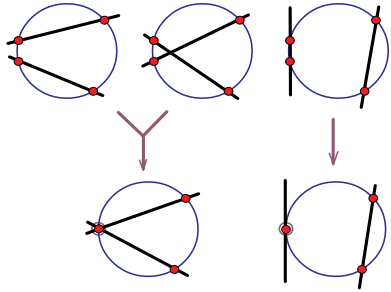}} 
\bigskip
\caption{Two distinct ways in which  parcellings degenerate.}
\end{figure}

Each time $L \in \mathcal Fact(d)$ produces in $\mathbb CP^2_\ast$ a union of lines with at least one pair of lines sharing its intersection point with the curve $\mathcal Q$, the point $\Pi(L) \in \mathcal D(\Pi_Q) =  \mathcal D(\Phi_Q)$. Moreover, as Figure 1 demonstrates, this change in cardinality of fibers occurs due to their bifurcations in the vicinity of $L$. As a result, points in $\mathcal D(\Phi_Q)$ give rise to effective divisors of degree $2d$ on the curve $\mathcal Q$ with at least one point in their support being of multiplicity greater than one. On the other hand,  any such divisor can be generated by intersecting a weighted set of lines with $\mathcal Q$ (just use any generalized parcelling). Note that any line tangent to $\mathcal Q$ also contributes a point of multiplicity two. However, the lines  tangent  to $\mathcal Q$ are not contributing to the bifurcation of the $\Pi_Q$-fibers (see the right diagram in Figure 1). Therefore, in order to conclude that any  divisor on $\mathcal Q$ with multiple points corresponds  to a point of $ \mathcal D(\Phi_Q)$, we need to use generalized parcellings which favor pairs of  lines that share their intersection with $\mathcal Q$ to a single tangent line (both patterns produce an intersection point of multiplicity two as shown in Figure 1). Evidently, this can be done, provided $d > 1$.

The requirement that $L \in \mathcal Fact(d)$ has a pair of linear forms vanishing at a point of 
$\{Q = 0\}$, locally, is  a single algebraic condition imposed on the coefficients of of the two forms.  Therefore, taking closures, it picks a codimension one subvariety  $\Pi_Q^{-1}(\mathcal D(\Pi_Q)) \subset Fact(d)$. Thus, $\mathcal D(\Pi_Q) \subset V_Q^\perp(d)^\circ$ is a subvariety of codimension  one as well.
Fortunately, since $\mathcal Q$ admits a rational parameterization by a  map 
$\alpha: \mathbb CP^1 \rightarrow \mathcal Q$,  the set  $\Pi_Q^{-1}(\mathcal D(\Pi_Q))$ can be described in terms of  solvability of a system of two rather simple equations. In homogeneous coordinates $[u_0 : u_1]$ such a parameterization $\alpha$ can be given  by the formula 
$$\alpha([u_0 : u_1]) = [\alpha_0([u_0 : u_1]) : \alpha_1([u_0 : u_1]): \alpha_2([u_0 : u_1])],$$ where $\alpha_0, \alpha_1, \alpha_2$ are some quadratic forms. For example, for $Q = x^2 + y^2 + z^2$, 
$$\alpha_0 = \mathbf{i}(u_0^2 - u_1^2),\; \alpha_1 = 2\mathbf{i}\, u_0 u_1, \; \alpha_2 = u_0^2 + u_1^2.$$ 
The inverse of $\alpha$ is produced by the central projection of $\mathcal Q$ onto a line in $\mathbb CP^2_\ast$ from a center located at $\mathcal Q$. Therefore, as \emph{abstract} algebraic curves,  
$\mathcal Q$ and $\mathbb CP^1$ are isomorphic. 

In fact, $L(x, y, z) \in \mathcal Fact(d)$   belongs to $\Pi_Q^{-1}(\mathcal D(\Pi_Q))$, if and only if, the system 
\begin{eqnarray}
\Big\{L(\alpha_0,\, \alpha_1,\, \alpha_2)\Big\}(u_0, u_1) = 0\nonumber\\
\Big\{\partial_x L(\alpha_0,\, \alpha_1,\, \alpha_2)\partial_{u_0} \alpha_0 + \partial_y L(\alpha_0,\, \alpha_1,\, \alpha_2)\partial_{u_0} \alpha_1 + \partial_z L(\alpha_0,\, \alpha_1,\, \alpha_2)\partial_{u_0} \alpha_2\Big\}(u_0 ,  u_1) = 0, 
\end{eqnarray}
that guarantees an existence of a multiple  zero for the polynomial $L(\alpha(u_0, u_1))$ in 
$\mathbb CP^1$, has a non-trivial solution $(u_0, u_1)$. Writing down explicitly the resultant of the two polynomials in the LHS of $(20)$, and thus the equation of $\Pi_Q^{-1}(\mathcal D(\Pi_Q))$, seems to be cumbersome. Note that the Euler identity 
$$2d\cdot L =  ( u_0 \partial_{u_0}\alpha_0  +   u_1\partial_{u_1}\alpha_0)\partial_x L + 
( u_0 \partial_{u_0}\alpha_1  +   u_1\partial_{u_1}\alpha_1) \partial_y L + 
( u_0 \partial_{u_0}\alpha_2  +   u_1\partial_{u_1}\alpha_2) \partial_z L$$
explains the asymmetry of $(20)$ with respect to the variable $u_0$: exchanging the roles of $u_0$ and $u_1$ leads to an equivalent system of equations.

Recall that for a projective variety $X$, the set  of zero-dimensional, degree $d$ effective divisors is the projective variety $Sym^d(X)$ (see  [Ch]). The map $\Phi_Q$ in formula $(16)$ induces a well-defined regular map of the varieties: 
\begin{eqnarray} 
\Psi_Q: Sym^d(\mathbb CP^2) \rightarrow  Sym^{2d}(\mathcal Q).
\end{eqnarray}
This map is produced by realizing a given multipole $w$, or rather its $\mathbb C^\ast$-orbit $\tilde w$, by a completely factorable polynomial $L(w)$ and then forming the intersection set $Z(L(w), Q)$ equipped with the appropriate multiplicities---the divisor $\Psi_Q(\tilde w)$. The map $\Psi_Q$ is onto: any effective divisor of degree $2d$ on $\mathcal Q$ admits a generalized parcelling, and thus is generated by intersecting $\mathcal Q$ with a weighted collection of lines. Since the number of generalized parcellings of $Z(L(w), Q)$ is finite, the map $\Psi_Q$ has finite fibers.
All this can be seen from a different angle. The divisor $\Psi_Q(\tilde w)$ uniquely determines the proportionality class of the function $L(w)|_{\{Q = 0\}}$---a point in $\mathbb CP(V_Q^\perp)$.  We have seen that the $\mathbb C^\ast$-equivariant map $\Phi_Q$ in $(16)$ is onto. Thus, 
$\Phi_Q/ \mathbb C^\ast: Sym^d(\mathbb CP^2) \rightarrow \mathbb CP(V_Q^\perp)$ is onto as well.  Since $\mathcal Q$ is a rational curve (topologically, a 2-sphere), there is a 1-to-1 correspondence between points of $\mathbb CP(V_Q^\perp) \approx Sym^{2d}(\mathbb CP^1_\ast)$ and of $Sym^{2d}(\mathcal Q)$. Moreover, as abstract varieties, $\mathbb CP(V_Q^\perp)$ and of $Sym^{2d}(\mathcal Q)$ are isomorphic. This isomorphism can be used to identify the maps $\Psi_Q$ and $\Phi_Q/ \mathbb C^\ast$. To simplify our notations, we will use the same symbol $\Psi_Q$ for both maps; when there is a need to distinguish them, we will just indicate the relevant target space.

Our combinatorial considerations imply that the locus $\mathcal D(\Psi_Q)$ consists of
the effective divisors of degree $2d$ on $\mathcal Q$ with at least one point in their support being of multiplicity at least two. While $\Psi_Q^{-1}(\mathcal D(\Psi_Q))$ is also described by $(20)$,  the locus $sing(\mathcal Fact(d))$ is the preimage of  effective divisors of degree $2d$ on $\mathcal Q$ with at least a pair of  points in their support being of multiplicity at least two (they can generate a double line), or at least one point being of multiplicity at least four (it corresponds to a double line tangent to 
$\mathcal Q$). 
  
Employing Lemma 4, we have established
\begin{lem} The ramification set $\mathcal D(\Phi_Q)$ for the map 
$\Phi_Q: \mathcal M(d) \rightarrow V_Q^\perp(d)$
is  the set $\{\Pi_Q(L)\}$ of complex codimension one, where $L \in \mathcal Fact(d)$ defines on  $\mathcal Q$ an effective divisor with  at least one point in its support being of  multiplicity  at least two. In other words, $\Pi_Q^{-1}(\mathcal D(\Phi_Q))$ is defined by the resultant of the two polynomials in the LHS of $(20)$. The ramification set  $\mathcal D(\Phi_Q)$ contains the $\Phi_Q$-image of the singular set $sing(\mathcal M(d))$. This image is of codimension two in $V_Q^\perp(d)$. It can be  identified with  completely factorable homogeneous  polynomials of degree $d$ on the the surface $\{Q = 0\}$ that have at least one pair of proportional linear factors. $\Box$
\end{lem}

Note that the space $Z$ of \emph{simple} effective degree $k$ divisors on the curve $\mathcal Q$ is homeomorphic to the space of simple effective degree $k$ divisors on the sphere $S^2 = \mathbb CP^1_\ast$. Therefore, $\pi_1(Z) \approx \mathsf B_k$, where $\mathsf B_k$ stands for the braid group in $k$ strings residing in the spherical shell  $S^2 \times [0, 1]$. In particular, $Sym^{2d}(\mathcal Q) \setminus \mathcal D(\Psi_Q)$ is a $K(\mathsf B_{2d}, 1)$-space.

Let $\mathsf B_{2d}^\#$ be  the preimage of the subgroup $S_{2d}^\# \subset S_{2d}$ in $(19)$ (of order $2^d \cdot d!$) under the canonical epimorphism $\mathsf B_{2d} \rightarrow S_{2d}$. We call it \emph{the braid group in} $2d$ \emph{strings with coupling}\footnote{Note that $\mathsf B_{2d}^\#$ is not  the  braid group in $2d$ strings colored with $d$ colors, each color marking a pair of strings!}.  It is a subgroup of index $(2d -1)!!$ in the braid group $\mathsf B_{2d}$.

Recall that a covering of a $K(\pi, 1)$-space is again a $K(\pi', 1)$-space, where $\pi'$ is an appropriate subgroup of $\pi$.  Since all our constructions are $\mathbb C^\ast$-equivariant, Lemmas 4 and 5 together with the arguments above lead to 

\begin{theorem} 
\begin{itemize}
\item The map 
$\Psi_Q: Sym^d(\mathbb CP^2) \rightarrow \mathbb CP^{2d} \approx 
V^\perp_Q(d)^\circ/ \mathbb C^\ast$ is a finite ramified covering with a generic fiber of cardinality 
$(2d - 1)!!$.  It is ramified over the subvariety 
$\mathcal D(\Psi_Q)$ whose points are the proportionality classes of   homogeneous polynomials  of degree $d$ on the quadratic surface $\{Q(x, y, z) = 0\}$ that define there effective divisors of degree $d$ with at least one multiple line\footnote{equivalently, that define on the curve $\mathcal Q \subset \mathbb CP^2_\ast$ effective divisors with a point of multiplicity at least two.}.

\item The space $\mathbb CP^{2d} \setminus \mathcal D(\Psi_Q)$ is a $K(\pi, 1)$-space with the group $\pi$ being isomorphic to the  braid group $\mathsf B_{2d}$ in $2d$ strings in the spherical shell $S^2 \times [0, 1]$. 

\item As a result, $Sym^d(\mathbb CP^2) \setminus \Psi_Q^{-1}(\mathcal D(\Psi_Q))$ is a  
$K(\mathsf B_{2d}^\#, 1)$-space, where  $\mathsf B_{2d}^\#$ is the braid group in $2d$ strings with coupling.

\item The  space $\mathcal M(d) \setminus \Phi_Q^{-1}(\mathcal D(\Phi_Q))$ also is  $K(\pi, 1)$-space with the group $\pi$ being isomorphic to an extension of the infinite cyclic group $\mathbb Z$ by the group $\mathsf B_{2d}^\#$.  $\Box$
\end{itemize}
\end{theorem}

Now, we will investigate the ramification locus of $\Phi_Q$ from a more refined point of view characteristic to the singularity theory.   First, we would like to understand  when a non-zero vector $w$ from the tangent cone $\mathcal T_L$ of 
$\mathcal Fact(d)$ at a point $L = \prod_{j =1}^d L_j$ is parallel to the subspace $V_Q(d)$, in other words, when $\mathcal T_L$ contains a vector $w$ that is mapped to zero under the projection
$\pi: V(d) \rightarrow V(d)/V_Q(d) \approx  V_Q^\perp(d)$. Away from the singularity set 
$sing(\mathcal Fact(d)) \subset \mathcal Fact(d)$, the $\Pi_Q$-image of such an $L$  belongs to a locus $\mathcal E \subset V_Q^\perp(d)$ over which 
$rank(\mathrm {d} \pi\big |_{\mathcal Fact(d)}) < dim(V_Q^\perp(d))$ at some point in 
$\Pi_Q^{-1}(\mathcal E)$. 
By the implicit function theorem, the ramification locus $\mathcal D(\Pi_Q) \subset V_Q^\perp(d)$ for the map $\Pi_Q : \mathcal Fact(d) \rightarrow   V_Q^\perp(d)$ (equivalently, for the map $\Phi_Q$) is contained in the union $\mathcal E \cup \Pi_Q(sing(\mathcal Fact(d)))$.  It can happen that at a singularity 
$L \in sing(\mathcal Fact(d))$ the tangent cone does not have vectors $w \neq 0$ with the property $\pi(w) = 0$, and still $\pi$ is ramified in the vicinity of  $\pi(L)$. For instance, consider the obvious projection $\pi$ of the real cone $x^2 + y^2 -z^2 = 0$ onto the $xy$-plane: $\pi$ is ramified at the origin $(0, 0)$, but for any $w \neq 0$ from the tangent cone at $(0, 0, 0)$, 
$\pi(w) \neq 0$.

Any $w \in \mathcal T_L$ is of the form 
$lim_{t \rightarrow 0} \frac{\prod_{j =1}^d (L_j + t M_j) - \prod_{j =1}^d L_j}{t}$, where each $M_j$ is 
an appropriate linear form in $x_j, y_j$, and $z_j$ or $M_j = 0 $ identically. This limit is equal to the polynomial $P = \sum_j (\prod_{i \neq j} L_i)M_j$. The vector $P$ at $L$ is parallel to the subspace $V_Q$ if and only if the polynomial $P$ is divisible by $Q$. In other words, the vector $P$ at $L$ is parallel to  $V_Q$ if and only if the polynomial $P$, being restricted to the curve $\mathcal Q$, is identically zero. Thus, we are looking for the $M_j$'s subject to the constraint: the polynomial 
$\sum_j (\prod_{i \neq j} L_i)M_j = (\prod_{i} L_i)(\sum_j  \frac{M_j}{L_j})$, being restricted to 
$\mathcal Q$,  is identically zero. For $L \neq 0$, this is equivalent to the constraint $(22)$ imposed 
on the rational functions $\{M_j/ L_j\}$:
\begin{eqnarray}
\sum_j  \frac{M_j}{L_j}\Big |_{\mathcal Q} = 0. 
\end{eqnarray}
Equation $(22)$ always have \emph{obvious} solutions: $\{M_j = \alpha_j L_j\}$, where 
$\{\alpha_j \in \mathbb C\}$ and $\sum_j \alpha_j = 0$. These are exactly solutions that represent the  zero tangent vector (the tip $L$ of the tangent cone). Indeed, put 
$\prod_{j =1}^d (L_j + t M_j) = \prod_{j =1}^d (L_j + t \alpha_j L_j)$ to conclude that $w = 0$ if and only if 
$\sum_j \alpha_j = 0$.

So, the proper question is how to describe all the $L \in \mathcal Fact(d)$ for which $(22)$ has a solution \emph{distinct} from the set of obvious solutions $\{M_j = \alpha_j L_j\}$ with  $\sum_j \alpha_j = 0$. The images of such $L$'s under the projection $\pi: V(d) \rightarrow  V_Q^\perp(d)$ will generate the locus 
$\mathcal E \subset V_Q^\perp(d)$  over which the differential  $\rm d\pi\big |_{\mathcal Fact(d)}$ is not of the maximal rank $dim(V_Q^\perp(d)) = 2d + 1$.  Evaluating the LHS of equation $(22)$ at 
$2d + 1$ generic points residing in $\mathcal Q$ imposes linear constraints on $3d$ variables---the coefficients of the $M_j$'s. If these constraints are independent, the solution space of the linear system 
must be of dimension $3d - (2d + 1) = d - 1$, which is exactly the dimension of the space formed by the obvious solutions. This indicates that, for a generic $L$, we should not expect any non-obvious solutions. Lemma 6 below validates  this guess.

Let us denote by $\mathbf M_L$ the quotient of the vector space of all solutions $\{M_j\}$ of $(22)$ by the subspace of obvious solutions, as defined above. 
The correspondence $L \Rightarrow dim(\mathbf M_L)$ gives rise to a new natural stratification of the space $\mathcal Fact(d)$, and thus, of the space $\mathcal M_Q(d)$.
In the new notations, $\pi(L) \in \mathcal E$ if and only if $\mathbf M_L \neq 0$. We suspect that  this stratification is consistent with, but cruder than the stratification induced by the orbit-types of the 
$\Sigma_d$-action.
\begin{lem} The loci $\mathcal E$ and $\mathcal D(\Pi_Q)$ coincide. As a result, the existence of a non-trivial solution for the system $(20)$ is equivalent to the existence of a non-obvious solution for the equation $(22)$.  Also, the locus $\mathcal E \supset \Phi_Q(sing(\mathcal M(d)))$.
\end{lem} 
{\bf Proof. \quad} We notice that if two distinct lines, say $\mathcal L_1$ and $\mathcal L_2$, share a point $p \in \mathcal Q$, then $(22)$ has a non-obvious solution $(M_1, M_2, 0, \; ... \;, 0)$. Indeed, 
inscribe in the quadratic curve $\mathcal Q$ any  "quadrilateral"  formed by the pair of lines 
$\mathcal L_1, \mathcal L_2$ together with a new pair of lines $\mathcal M_1, \mathcal M_2$. The lines 
$\mathcal L_1$, $\mathcal M_1$ share a point $q \in \mathcal Q$,  the lines $\mathcal L_2$, 
$\mathcal M_2$ share a point $r \in \mathcal Q$, and the lines $\mathcal M_1$, 
$\mathcal M_2$ share a point $s\in \mathcal Q$, all four points $p, q, r, s$ being distinct. Next, pick some linear forms $M_1$ and $M_2$ representing $\mathcal M_1$ and $\mathcal M_2$. Then one can find constants $\lambda_1, \lambda_2$ so that  the polynomial 
$Q = \lambda_1 M_1 L_2 + \lambda_2 M_2 L_1$. 
The argument is very similar to the one used in the proof of  Lemma 1. Thus, 
$\lambda_1 M_1 L_2 + \lambda_2 M_2 L_1\big|_\mathcal Q = 0$, which can be written in the form 
$\lambda_1 M_1/ L_1 + \lambda_2 M_2 /L_2\big|_\mathcal Q = 0$ required by $(22)$. Here, evidently, 
$M_1$ is not proportional to $L_1$ and $M_2$ is not proportional to $L_2$. 
Therefore, for any $L \in \mathcal Fact(d)$ that defines on $\mathcal Q$ an effective divisor containing a point of multiplicity two,  $\Pi_Q(L)$ must  belong to the locus $\mathcal E$. It follows "by continuity" that all the $L$'s that produce multiplicity functions $\mu: Z(L, Q) \rightarrow \mathbb N$, distinct from the identity function $1$, project to $\mathcal E$ via $\Pi$. According to Lemma 5, these are exactly the factorable polynomials that project to the "combinatorial'" ramification locus $\mathcal D$. $\Box$

It follows from Lemma 6 that, for a generic $P \in V(d)$, the affine subspace 
$P + V_Q(d)$  hits transversally the subvariety ${\mathcal Fact}(d)$ at a finite set of points whose cardinality is exactly $\kappa(d)$. Therefore, $\kappa(d)$ must be the degree of that variety.

Let $X, Y$ be quasi-affine varieties and let $F: X \rightarrow Y$ be a proper regular map.
In Theorem 7 below, the ramification set $\mathcal D(F)$ of a mapping $F: X \rightarrow Y$ with finite fibers is understood to be the closure of a set $\mathcal D^\circ(F)$. By definition, $y \in \mathcal D^\circ(F)$ when $F^{-1}(y)$ contains a point $x$ such that  there  is a non-zero vector 
$v \in \mathcal T_xX$ that is mapped to zero in $\mathcal T_yY$ by the differential $DF$. 

The arguments above lead  to our main result:
\begin{theorem} 
\begin{itemize} 
\item The map $ \Phi_Q:  \mathcal M(d) \rightarrow {\mathcal Fact}_Q(d) = V_Q^{\perp}(d)^\circ$ is onto, and its generic fiber  is a finite set of cardinality $(2d - 1)!!$. This map  takes the multipole space $\mathcal M(d)$---the space of a 
 $\mathbb C^\ast$-bundle $\eta (k)$ over the variety $Sym^d(\mathbb CP^2)$---onto the complex vector space $V_Q^{\perp}(d)$ of dimension $2d + 1$ with its origin being deleted.
The map $\Phi_Q$ extends to a map 
$\overline \Phi_Q : \overline \mathcal M(d) \longrightarrow  V_Q^{\perp}(d)$ with finite fibers.
 
\item $\Phi_Q$ is ramified over the discriminant variety $\mathcal D( \Phi_Q)$---the set of polynomials 
$P \in {\mathcal Fact}_Q(d)$  whose zero sets are effective divisors on the surface 
$\{Q(x, y, z) = 0\}$ with at least one of their line components being of multiplicity at least
two.\footnote{Equivalently, the set of polynomials $P$ for which $\mathcal P \cap \mathcal Q$ has a point of multiplicity at least two.}
The subvariety $\mathcal D( \Phi_Q) \subset V_Q^{\perp}(d)^\circ$ is of  complex codimension one,  and is described in terms of solvability of $(20)$ or $(22)$.  It contains the 
$\Phi_Q$-image of the singular set $sing(\mathcal M(d))$ as a codimension one subvariety.

\item In fact,  $\Phi_Q$, $\overline \Phi_Q$ are  $\mathbb C^\ast$-equivariant maps. Therefore, $\Phi_Q$ gives rise to a surjective map 
$$\Psi_Q : Sym^d(\mathbb CP^2) \rightarrow \mathbb CP^{2d} = \mathbb CP(V_Q^{\perp}(d)^\circ)$$ 
of degree $(2d -1)!!$ which  is  also ramified over a subvariety 
$\mathcal D(\Psi_Q) \subset \mathbb CP^{2d}$ of codimension one. The compliment 
$\mathbb CP^{2d} \setminus \mathcal D(\Psi_Q)$ is a $K(\mathsf B_{2d}, 1)$-space, 
$Sym^d(\mathbb CP^2) \setminus \Psi_Q^{-1}(\mathcal D(\Psi_Q))$ is a $K(\mathsf B_{2d}^\#, 1)$-space,  where $\mathsf B_{2d}^\#$ is the braid group with coupling, while 
$\mathcal M(d) \setminus \Phi_Q^{-1}(\mathcal D(\Phi_Q))$ is a $K(\pi, 1)$-space,  where $\pi$ an extension of  $\mathbb Z$ by  $\mathsf B_{2d}^\#$.

\item The degree of the variety  of completely factorizable homogeneous polynomials ${\mathcal Fact}(d) \subset V(d)$ is also $(2d - 1)!!$. This variety is invariant under  the obvious $\mathbb C^\ast$-action on $V(d)$. \qquad $\Box$
\end{itemize}
\end{theorem}
Theorem 7 has a number of topological implications, some of them dealing with interesting ramifications over complex projective spaces. Our next goal is to describe these implications.

Let the space $X$ be of a homotopy type of a connected, finite-dimensional CW-complex. Recall that the Dold-Thom Theorem [DT] 
links the homotopy groups of $Sym^d(X)$, $d$ being large, with the integral homology of a space $X$.
By picking a base point $a$ in $X$, one gets a stabilization map $Sym^d(X) \rightarrow Sym^{d + 1}(X)$, well-defined by the formula $(x_1, x_2, ... , x_d) \rightarrow (a, x_1, x_2, ... , x_d)$. Here $\{x_i \in X\}$. This provides us with canonical homomorphisms $\pi_k(Sym^d(X)) \rightarrow \pi_k(Sym^{d + 1}(X))$ of the $k$-th homotopy groups.

In our case, the Dold-Thom Theorem claims that $lim_{d \rightarrow \infty}\; \pi_k(Sym^d(\mathbb CP^2)) = H_k(\mathbb CP^2; \mathbb Z)$. 
In particular, $lim_{d \rightarrow \infty}\; \pi_2(Sym^d(\mathbb CP^2)) = \mathbb Z = lim_{d \rightarrow \infty}\; \pi_4(Sym^d(\mathbb CP^2))$. Also,  $lim_{d \rightarrow \infty}\pi_k(Sym^d(\mathbb CP^2)) = 0$,
provided $k = 1, 3$ or $k > 4$.  On the other hand, $\pi_2(\mathbb CP^{2d}) = \mathbb Z$, but 
$\pi_4(\mathbb CP^{2d}) = 0$. Therefore, at least for large $d$'s, there is an infinite order element 
$ \alpha \in \pi_4(Sym^d(\mathbb CP^2))$ that is mapped by $(\Psi_Q)_\ast$ to zero  
and an element $\beta \in \pi_2(Sym^d(\mathbb CP^2))$, so that  $(\Psi_Q)_\ast(\beta) \in\pi_2(\mathbb CP^{2d})$ is a generator. It is possible to realize $\alpha$ and $\beta$ geometrically. What is clear ratherway, that the spheroids  $\alpha$ and $\beta$ can not be pushed into the aspherical portion $Sym^d(\mathbb CP^2) \setminus \Psi_Q^{-1}(\mathcal D(\Psi_Q))$ of $Sym^d(\mathbb CP^2)$. 
The realization of $\alpha$ is based on an interesting fact that I learned from Blaine Lawson: the quotient of $\mathbb CP^2$ by the complex conjugation 
$\tau : [x : y : z] \rightarrow [\bar x : \bar y : \bar z]$ is homeomorphic to the sphere $S^4$. Therefore, the map $\phi: \mathbb CP^2 \rightarrow \mathbb CP^2 \times \mathbb CP^2$ given by the formula 
$\phi (p) = (p, \tau(p))$, where $p \in \mathbb CP^2$, is evidently $\mathbb Z_2$-equivariant with respect to the $\tau$-action in the domain and the symmetrizing action in the range. This gives rise to the desired quotient map $\alpha: S^4 \approx \mathbb CP^2 / \{\tau\} \rightarrow  Sym^2(\mathbb CP^2)$ that survives into the higher symmetric powers of $\mathbb CP^2$. Constructing  class $\beta$ is straightforward: it is given by the obvious inclusion $S^2 \approx \mathbb CP^1 \subset \mathbb CP^2$ followed by the diagonal embedding $\Delta: \mathbb CP^2 \rightarrow Sym^d(\mathbb CP^2)$.

The existence of a non-trivial $\alpha : S^4 \rightarrow Sym^d(\mathbb CP^2)$ whose $\Psi_Q$-image is null-homotopic in $ \mathbb CP^{2d}$ has a curious implication:
\begin{cor} 
For any $d \geq 2$, there is a family of  polynomial  functions of degree $d$ on the quadratic surface $\{Q = 0\}$ that is parameterized by a $5$-dimensional disk which does not admit a continuous lifting to the multipole space $\mathcal M(d)$. At the same time, the functions parameterized  by the $4$-sphere forming the boundary of the disk can be continuously represented by multipoles. $\Box$
\end{cor}
Any map $f: X \rightarrow Y$ induces a natural map $f_{\ast}^k: Sym^k(X) \rightarrow Sym^k(Y)$ that is defined by the formula  $(f_{\ast}^k)[\sum_\nu \mu_\nu x_\nu] = \sum_\nu \mu_\nu f(x_\nu)$, where $\{x_\nu \in X\}$ and $\{\mu_\nu \in \mathbb N\}$. When $f$ is $1$-to-$1$ or onto, so is $f_{\ast}^k$.

The construction of $f_{\ast}^k$ provides us with a rich source of interesting ramified coverings. Consider for example, a semi-free  cyclic action  on a sphere $S^2 \subset \mathbb R^3$.  The group $\mathbb Z_l$ acts on $S^2$ by rotations  around a fixed axis on the angles  that are multiples of $2\pi/l$. Topologically, the orbit-space $\tilde S^2 : = S^2 / \mathbb Z_l$ is again a 2-sphere. Let $f: S^2 \rightarrow \tilde S^2$ be the orbit-map. Then $f_{\ast}^k :  Sym^k(S^2) \rightarrow Sym^k(\tilde S^2)$ gives an example of a ramified degree $l^k$ covering map of $\mathbb CP^k$ over itself!  For $l = 2$, we shall see later how this degree $2^k$ ramification $f_{\ast}^k :   \mathbb CP^k  \rightarrow  \mathbb CP^k$ is linked to the multipole spaces on quadratic surfaces. On the other hand, by a simple cohomological  argument, any map $F: \mathbb CP^k \rightarrow \mathbb CP^k$ has a degree that is the $k$-th power of a non-negative integer $l$.  For instance, there is no ramified map from $\mathbb CP^2$ to itself of degrees that are not of the form $l^2$.

{\bf Question } Given a  closed oriented manifold $M$,  what are  possible degrees  of  maps from $M$ to itself? Evidently, the answer depends on the homotopy type of $M$.

When $\mathbb Z_2$ acts freely on $S^2$ by the central symmetry, the orbit-map $f: S^2 \rightarrow \mathbb RP^2$ gives rise to another interesting ramified map $f_{\ast}^k :  \mathbb CP^k \rightarrow \mathbb RP^{2k}$ of degree $2^k$ (its existence follows from our results in Section 3 dealing with real multipole spaces).

In order to derive next few corollaries of Theorem 7, we need to take a  detour aimed at constructing (with the help of $f$) transfer maps that take effective $0$-divisors on $Y$ to effective $0$-divisors on $X$. These constructions are variations on the theme of the classical Hurwitz' Theorem (see [H], pp. 299-304). Our next goal is to present constructions and arguments that lead to Theorem 10.\smallskip

Let  $X$ and $Y$ be smooth complex projective varieties and $f: X \rightarrow Y$ a regular surjective  map with finite fibers. 
With any $x \in X$ we associate a multiplicity number $\mu_f(x)$. It is the multiplicity attached to the intersection of the $f$-graph $\Gamma_f $ with the subspace $X \times f(x) \subset X \times Y$ at the point $(x, f(x))$. Since  all the $f$-fibers are finite, the intersection $\Gamma_f \cap (X \times f(x))$ is finite as well. 

Next, with any $y \in Y$ we associate an effective $0$-divisor $D_{f^{-1}(y)} := \sum_{x \in f^{-1}(y)} \mu_f(x) x$ whose degree is $y$-independent and coincides with the degree $d_f$ of the map $f$. The correspondence $y \Rightarrow D_{f^{-1}(y)}$ produces a regular embedding 
\begin{eqnarray}
f^\#: Y \rightarrow Sym^{d_f}(X)
\end{eqnarray}
For any $k$,  the map $f^\#$  gives rise  to an embedding
\begin{eqnarray}
f^\#_k: Sym^k(Y) \rightarrow Sym^{k\cdot d_f}(X)
\end{eqnarray}
defined by the formula $f^\#_k\big(\sum_\nu \mu(y_\nu)y\big) = \sum_\nu \mu(y_\nu)D_{f^{-1}(y)}$.

Therefore, applying this construction to the setting of Theorem 7 with $X = Sym^d(\mathbb CP^2)$,  $Y = \mathbb CP^{2d}$, and $f = \Psi_Q$, we get the following proposition:
\begin{cor} 
Any irreducible quadratic form $Q(x, y, z)$ determines  canonical  embeddings 
$$\Psi^\#_{Q, k} : Sym^k(\mathbb CP^{2d}) \rightarrow Sym^{k \cdot (2d - 1)!!}(Sym^d(\mathbb CP^2)),$$
where $d$ and $k$ are arbitrary whole numbers.
In particular, $\Psi_Q^\# := \Psi^\#_{Q, 1}$ embeds  $\mathbb CP^{2d}$ into \hfill\break
$Sym^{(2d - 1)!!}(Sym^d(\mathbb CP^2)).$ \qquad $\Box$
\end{cor}
For a map $f: X \rightarrow Y$ as above, one can define a  map $f_\#$ that takes effective $0$-divisors on $X$ into effective $0$-divisors on $X$.  By definition, each point $x \in X$ is mapped to the divisor $D_{f^{-1}(f(x))}$. By linearity, we have 
\begin{eqnarray}
f_\# \Big(\sum_\nu \mu(x_\nu) x_\nu\Big) = \sum_\nu \mu(x_\nu) D_{f^{-1}(f(x_\nu))}.
\end{eqnarray}
This map transforms  $0$-divisors of degree $k$  into $0$-divisors of degree $k\cdot d_f$: 
\begin{eqnarray}
f_\#^k: Sym^k(X) \rightarrow Sym^{k\cdot d_f}(X)
\end{eqnarray}
Both maps, $f_\#^k$ from $(26)$ and $f^\#_k$ from $(24)$, have the same targets. Moreover, their images coincide. Indeed, for each point $y \in Y$, pick any point $x \in f^{-1}(y)$. Then,  
$f^\#(y) = D_{f^{-1}(y)} = f_{\#}(x)$. Recall, that $f^\#_k$ is a 1-to-1 map and thus is invertible over its image. Therefore, the  map
\begin{eqnarray}
(f^\#_k)^{-1} \circ f_\# ^k: Sym^k(X) \rightarrow Sym^k(Y)
\end{eqnarray}
is well-defined.
\begin{lem}  Let $X$ and $Y$ be smooth complex projective varieties and $f: X \rightarrow Y$ be a regular onto map with finite fibers. The  map $(f^\#_k)^{-1} \circ f_\# ^k$ in $(27)$ coincides with the natural map $f_{\ast}^k: Sym^k(X) \rightarrow Sym^k(Y)$. It defines a ramified covering with a generic fiber of cardinality $(d_f)^k$.
\end{lem} 
{\bf Proof \quad} The map $f_\#$ takes each point $x \in X$ to the divisor 
$D_{f^{-1}(f(x))}$.\footnote{Recall  that the multiplicity of $x$ in $D_{f^{-1}(f(x))}$ is $\mu_f(x)$.}. At the same time, the transfer $f^\#$ takes $f(x)$ to the same divisor $D_{f^{-1}(f(x))}$. Thus, $(f^\#)^{-1}\circ f_\#$ maps $x$ to $f(x)$. Extending this argument by linearity proves the claim. $\Box$

Now we examine how these constructions apply to projective curves in $\mathbb CP^2_\ast$ and eventually to the multipole spaces. The role of $X$ will be played  by a curve $\mathcal C \subset \mathbb CP^2_\ast$ (most importantly, by $\mathcal Q$), the role of $Y$ by a pencil of lines in 
$\mathbb CP^2_\ast$ through a point $p$.  The map  $f$ takes any point $q \in \mathcal C$ to a line $\mathcal L$ passing through $q$ and $p$.

Any linear  embedding  $\rho: \mathbb CP^1 \subset \mathbb CP^2$ induces an embedding $\rho^d : Sym^d(\mathbb CP^1) \rightarrow Sym^d(\mathbb CP^2)$. The geometry of $\rho^d$ is tricky.  Since the quotient space $\mathbb CP^2 / \mathbb CP^1$ is homeomorphic to a 4-sphere $S^4$,  we get a homeomorphism $$Sym^d(\mathbb CP^2) / Sym^d(\mathbb CP^1) \approx Sym^d(S^4),$$ 
but even the spaces $\{Sym^d(S^4)\}$ have subtle topology. For example, $Sym^2(S^4)$ is homeomorphic to a mapping cylinder of a map $\Sigma^4(\mathbb RP^3) \rightarrow S^4$, where $\Sigma^4(\sim)$ denotes the fourth suspension (see [Ha], Example 4K.5).

The regular   map
\begin{eqnarray} 
\hat\Psi_Q: Sym^d(\mathbb CP^1) \stackrel{\rho^d}{\rightarrow} Sym^d(\mathbb CP^2)  \stackrel{\Psi_Q}{\rightarrow} Sym^{2d}(\mathcal Q) \approx \mathbb CP^{2d} 
\end{eqnarray}
describes the role and place of $\mathbb C$-\emph{planar} multipoles (they are linked to the planarity of the quadra- poles and octapoles in the deconstruction of  CMBR that was briefly mentioned in the introduction).
In contrast with $\Psi_Q$, we will see that the map $\hat\Psi_Q$ is 1-to-1. Since  $Sym^d(\mathbb CP^1)$ is homeomorphic to $\mathbb CP^d$, its image in $Sym^{2d}(\mathcal Q)$ is homeomorphic to $\mathbb CP^{d}$ as well (compare this with $(21)$). Moreover, the $2d$-cycle 
$\hat\Psi_{Q\; \ast}[Sym^d(\mathbb CP^1)]$ is homologous to the cycle $[Sym^d(\mathcal Q)]$  (equivalently, to the cycle $[\mathbb CP^d]$) defined by a natural  imbedding $Sym^d(\mathcal Q) \subset Sym^{2d}(\mathcal Q)$ (equivalently, by a linear imbedding $\mathbb CP^d \subset \mathbb CP^{2d}$).

Recall that  points $L  \in  \mathbb CP^1 \subset \mathbb CP^2$ correspond to a pencil of lines $\mathcal L \subset \mathbb CP^2_\ast$ that pass through a particular point $p \in \mathbb CP^2_\ast$. That $p$ determines the embedding $\rho$. Each point $q \in \mathcal Q$ determines a unique line $\mathcal L_q$ that passes through $q$ and $p$, and therefore, a unique point $L_q$ in the dual subspace $\mathbb CP^1 \subset \mathbb CP^2$. Let $f: \mathcal Q \rightarrow \mathbb CP^1$ be a $2$-to-$1$ ramified map defined by the correspondence $q \Rightarrow L_q$. As described in (24) and (26), $f$ gives rise to maps $f_d^\# : Sym^d(\mathbb CP^1) \rightarrow Sym^{2d}(\mathcal Q)$ and $f^d_\# : Sym^d(\mathcal Q) \rightarrow Sym^{2d}(\mathcal Q)$ with  $f_d^\#$ being a $1$-to-$1$ map. In fact, examining the construction of $\Psi_Q$, we see that $\hat\Psi_Q = f_d^\#$. Moreover, using Lemma 7, the map $f_\ast^d : Sym^d(\mathcal Q) \rightarrow Sym^d(\mathbb CP^1)$ factors as $(f^\#_d)^{-1} \circ f_\# ^d$. 
Therefore, the ramification locus $\mathcal D(f_\ast^d) \subset Sym^d(\mathbb CP^1)$ for $f_\ast^d$ is the $\hat\Psi_Q^{-1}$-image of the ramification locus 
$\mathcal D(f_\#^d) \subset f_\#^d(Sym^d(\mathcal Q)) \subset Sym^{2d}(\mathcal Q)$ for $f_\#^d$. 

Figure 2 illustrates the case $d = 3$. Passing from the first to the second column depicts the map $f^3_\#: Sym^3(\mathcal Q) \rightarrow Sym^6(\mathcal Q)$ generated by the linear projection from a center $p$ located at infinity.  Passing from the first to the third column depicts the map $f^3_\ast: Sym^3(\mathcal Q) \rightarrow Sym^3(\mathbb CP^1) \approx \mathbb CP^3$. Here $\mathbb CP^1$ is viewed as the pencil of lines in $\mathbb CP^2_\ast$ through $p$. Topologically, the map $f^3_\ast$ is a 8-to-1 ramification of $\mathbb CP^3$ over itself. It is described in some detail in Example 1 and is depicted in Figure 3. The passage from the second to the third column is a 1-to-1 correspondence.
\begin{figure}[ht]
\centerline{\includegraphics[height=5in,width=4in]{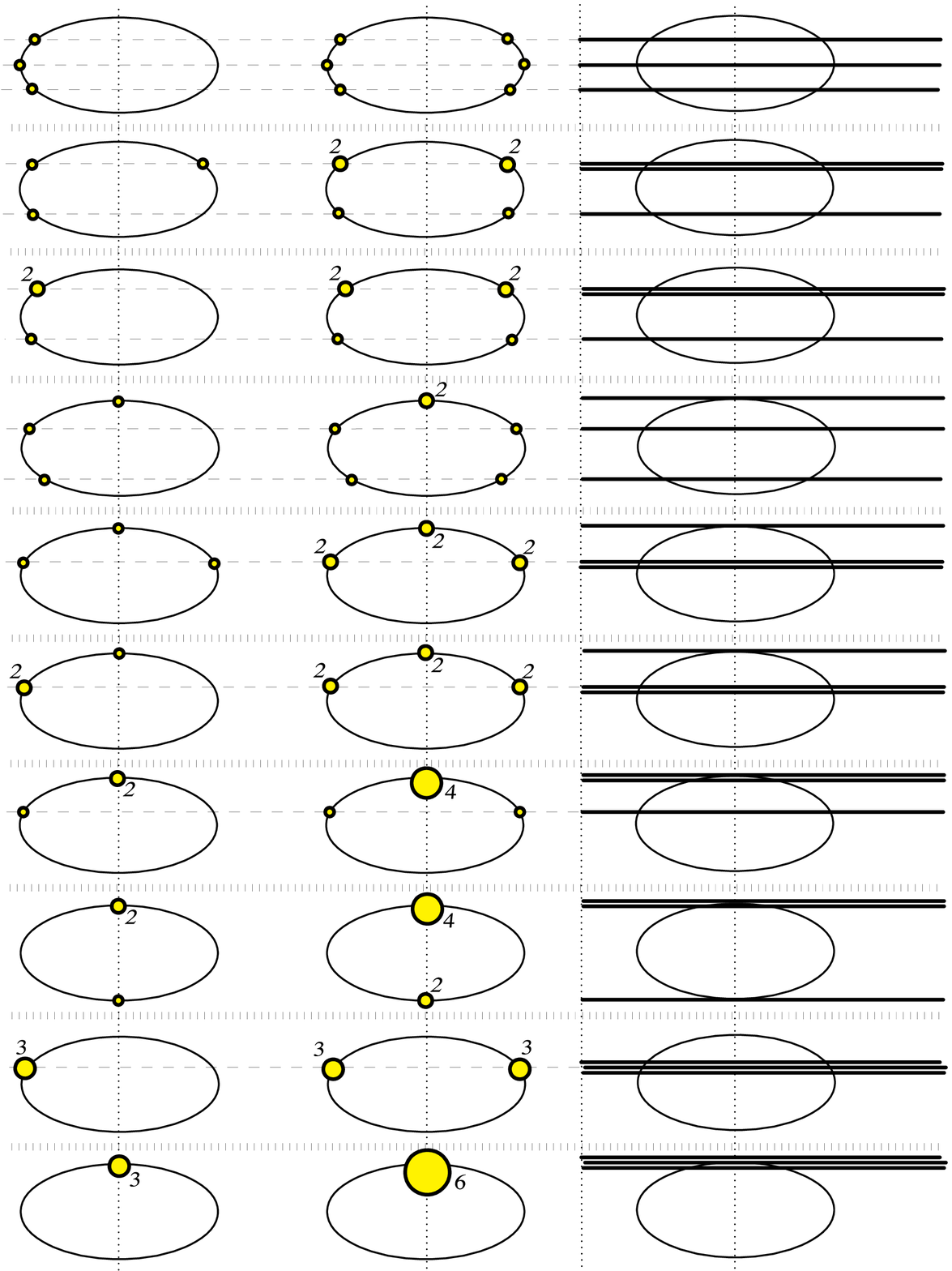}} 
\bigskip
\caption{}
\end{figure}

When $\mathcal L_q$ is not tangent to $\mathcal Q$ at $q$, then let  $q^\star \neq q$ be "the other point" in $\mathcal Q$ that belongs to the line $\mathcal L_q$.  When $\mathcal L_q$ is tangent to  $\mathcal Q$ at $q$, then by definition, put $q^\star = q$. The correspondence $\tau_p : q \Rightarrow q^\star$ is an involution on $\mathcal Q$ with two fixed points $a$ and $b$. Its orbit-space $\mathcal Q / \{\tau_p\}$ topologically is a  $2$-sphere. 
The involution $\tau_p$ induces an involution $\tau_p^k$ that acts on the space $Sym^k(\mathcal Q)$. By definition, any effective divisor $\sum_\nu \mu(q_\nu) q_\nu \in Sym^k(\mathcal Q)$ is transformed  by  $\tau_p^k$ into the divisor $\sum_\nu \mu(q_\nu) q_\nu^\star$.
Evidently, the image of $f^d_\# : Sym^d(\mathcal Q) \rightarrow Sym^{2d}(\mathcal Q)$ is contained in the $\tau_p^{2d}$-invariant part of the space $Sym^{2d}(\mathcal Q)$; however, not any invariant divisor belongs to that image. Each divisor from $Im(f^d_\#)$ not only must be $\tau_p^d$-invariant, but in addition, its multiplicity at the fixed points $a$ and $b$ must be \emph{even}. 
In other words, such a divisor $D^\#$ must be of the form $$2\mu_a a + 2\mu_b b + \sum_\nu \mu(q_\nu) [q_\nu +  q_\nu^\star],$$ where $q_\nu \neq a, b$ and $\mu_a, \mu_b \in \mathbb Z_+$. The cardinality of the $f^d_\#$-fiber over $D^\#$ is given by the formula 
\begin{eqnarray} 
|(f^d_\#)^{-1}(D^\#)| = \prod_\nu [\mu(q_\nu)  + 1]. 
\end{eqnarray}
Indeed, for any $q \in \mathcal Q \setminus \{a, b\}$ there are $\mu(q_\nu)  + 1$ effective divisors of degree $\mu(q_\nu)$ with the support in $q_\nu \coprod  q_\nu^\star$ and whose $f^{\mu(q_\nu)}_\#$-image is $\mu(q_\nu) [q_\nu +  q_\nu^\star]$; at the same time,  there is a unique divisor $\mu_a a$ whose $f^{\mu(\mu_a)}_\#$-image is $2\mu_a a$.
When all $\mu_\nu = 1$ and $\mu_a = 0 = \mu_b$, the fiber  $(f^d_\#)^{-1}(D^\#)$ is of cardinality $2^d$.
\smallskip

Now, we are in position to justify previous claims about  the image of the map $\hat\Psi_Q$ in (28).
Pick a divisor $D_0$ of degree $d$ in $\mathcal Q$ such that the supports of $D_0$ and $\tau_p^d(D_0)$ do not share common points. For any divisor $D$ of degree $d$ on $\mathcal Q$ the operation $D \Rightarrow D + D_0$ gives rise to an imbedding $I_{D_0}: Sym^d(\mathcal Q) \subset Sym^{2d}(\mathcal Q)$. We claim that  the intersection 
$\hat\Psi_Q(Sym^d(\mathbb CP^1)) \cap I_{D_0}(Sym^d(\mathcal Q))$ is a singleton. Indeed, $\hat\Psi_Q(Sym^d(\mathbb CP^1))$ consists of   $\tau_p^d$-invariant divisors, but the only invariant divisor of the form $D + D_0$ is the divisor $D_0 + \tau_p^d(D_0)$. Therefore, using the cohomology ring structure of $Sym^{2d}(\mathcal Q) \approx \mathbb CP^{2d}$, the cycle 
$[\hat\Psi_Q(Sym^d(\mathbb CP^1))] \in H_{2d}(Sym^{2d}(\mathcal Q); \mathbb Z)$ must be equal to the cycle   $[I_{D_0}(Sym^d(\mathcal Q))]$. Using the diffeomorphism $\mathcal Q \approx \mathbb CP^1$, this translates into the generating cycle 
$[\mathbb CP^d] \in H_{2d}(\mathbb CP^{2d}; \mathbb Z)$.
\smallskip

Examining $(29)$, we see that the ramification set $\mathcal D(f_\#^d)$ for 
$f_\#^d : Sym^d(\mathcal Q) \rightarrow Sym^{2d}(\mathcal Q)$ is comprised of divisors of two kinds: 1) the ones that contain at least one summand  of the form $2(q + q^\star)$, where $q \neq a, b$, and 2) the ones that   contain $2a$ or  $2b$ as a summand. 
For example, if $D \in Sym^d(\mathcal Q)$ contains a pair of  distinct points $q_1,  q_2 = q_1^\star$ and the rest of points  in the support of $D$ are generic (i.e., for $i, j > 2$,   $q_i \neq q_j^\star$), then $(f_\#^d)^{-1}(f_\#^d (D))$ consists of $2^d - 2^{d - 2} = 3\cdot 2^{d - 2}$ elements. 
At the same time, as we perturb $D$ in order  to avoid the coincidence $q_2 = q_1^\star$,\, $(f_\#^d)^{-1}(f_\#^d (D))$ consists of $2^d$ elements.

Therefore, in view of Lemma 7, the ramification set 
$\mathcal D(f_\ast^d)$ for $f_\ast^d : Sym^d(\mathcal Q) \rightarrow Sym^d(\mathbb CP^1)$ is comprised of divisors of two kinds: 1) the ones that contain at least one summand  of multiplicity at least two,  and 2) the ones that   contain points $L_a$ or  $L_b$ giving rise to lines $\mathcal L_a$ and  $\mathcal L_b$ passing through the point $p$ and tangent to the curve $\mathcal Q$. 

Given a space $X$, let us denote by $\Delta_d(X) \subset Sym^d(X)$ the discriminat set formed by the divisors containing points of multiplicity at least two. Also, for any point $a \in X$, we denote by 
$Sym^{d - k}_{ka}(X)$ the subset of $Sym^d(X)$ formed by the divisors containing the summand 
$k\cdot a$.

Thus, 
\begin{eqnarray} 
\mathcal D(f_\ast^d) = \Delta_d(\mathbb CP^1) \cup Sym^{d - 1}_a(\mathbb CP^1) \cup 
Sym^{d - 1}_b(\mathbb CP^1)
\end{eqnarray}
and 
 \begin{eqnarray} 
 \mathcal D(f_\#^d) = \Big\{\Delta_{2d}(\mathcal Q) \cup Sym^{2d -2}_{2a}(\mathcal Q) \cup 
 Sym^{2d -2}_{2b}(\mathcal Q)\Big\}^{\tau_p}.
 \end{eqnarray}
Employing  the Vi\'{e}te Map $V: Sym^d(\mathbb CP^1) \rightarrow \mathbb CP^d$,  we transplant the algebraic set $\mathcal D(f_\ast^d)$ into the space $\mathbb CP^d$. The image of $\Delta_d(\mathbb CP^1)$ under the Vi\'{e}te Map is the classical discriminant variety in $\mathcal D_d \subset \mathbb CP^d$, while the images of $Sym^{d - 1}_a(\mathbb CP^1)$ and of $Sym^{d - 1}_b(\mathbb CP^1)$ form  linear subspaces $\mathbb CP^{d - 1}_a$ and  
$\mathbb CP^{d - 1}_b$ in $\mathbb CP^d$. Thus, $V(\mathcal D(f_\ast^d)) = \mathcal D_d \cup \mathbb CP^{d - 1}_a \cup \mathbb CP^{d - 1}_b$. It follows from [K], Theorem 6.1, that each of the two spaces 
$\mathbb CP^{d - 1}_a$ and $\mathbb CP^{d - 1}_b$ is \emph{tangent} to the discriminant variety $\mathcal D_d$ along, respectively,  the  linear  subspaces $\mathbb CP^{d - 2}_a$ and 
$\mathbb CP^{d - 2}_b$. 

Note that the complement $Sym^d(\mathbb CP^1) \setminus \mathcal D(f_\ast^d)$ is the configuration space of $d$-tuples of distinct points  in the domain $\mathbb CP^1 \setminus (a^\ast \cup b^\ast)$. Here $a^\ast,   b^\ast$ stand for the two points in $\mathbb CP^1$ that are dual to the lines passing through $p$ and tangent to the curve $\mathcal Q$. Therefore, it is a $K(\mathsf B_d^{ann}, 1)$-space, where $\mathsf B_d^{ann}$ stands for the braid group in $d$ strings residing in a cylinder with an annulus base $[S^2 \setminus (a^\ast \cup b^\ast)] \times [0,1]$.

Using the birational identifications 
$Sym^d(\mathcal Q) \approx \mathbb CP^d \approx Sym^d(\mathbb CP^1)$, we have constructed a ramified covering $\Gamma_Q: \mathbb CP^d  \approx Sym^d(\mathcal Q) \stackrel{f^d_\ast}{\rightarrow} Sym^d(\mathbb CP^1) \approx \mathbb CP^d$. The following proposition summarizes the conclusions of our arguments above (centered on $(28)$ and $(30), (31)$). It describes an intricate stratified geometry of this ramified covering, a geometry that  is ``reductive" in its nature with respect to the shift $d \Rightarrow d - 1$.  
\begin{theorem} 
\begin{itemize}
\item Any irreducible quadratic form $Q$ gives rise to  a   ramified covering \hfill\break$\Gamma_Q: \mathbb CP^d \rightarrow \mathbb CP^d$ of degree $2^d$.  
The map $\Gamma_Q$ is  ramified over the algebraic set 
$$\mathcal D(\Gamma_Q) := \mathcal D_d \cup \mathbb CP^{d -1}_a \cup \mathbb CP^{d -1}_b$$--- the  Vi\'{e}te image of the set 
$\Delta_d(\mathbb CP^1) \cup Sym^{d - 1}_a(\mathbb CP^1) \cup  Sym^{d - 1}_b(\mathbb CP^1).$ 

\item The discriminant variety $\mathcal D_d \subset  \mathbb CP^d$ is $Q$-independent. The two linear subspaces $\mathbb CP^{d - 1}_a$ and $\mathbb CP^{d - 1}_b$ are tangent to the  variety $\mathcal D_d$ along, respectively, subspaces $\mathbb CP^{d - 2}_a$ and $\mathbb CP^{d - 2}_b$. 
A generic point of $\mathcal D_d$ has a $\Gamma_Q$-fiber of cardinality 
$3\cdot 2^{d - 2}$,\footnote{that is, $2^{d - 2}$ less than a generic $\Gamma_Q$-fiber.} while a generic  point of $\mathbb CP^{d -1}_a \cup \mathbb CP^{d -1}_b$ has a fiber of cardinality $2^{d - 1}$.

\item The complement $\mathbb CP^d \setminus \mathcal D(\Gamma_Q)$ to the ramification set is a $K(\mathsf B_d^{ann}, 1)$-space, where $\mathsf B_d^{ann}$ denotes the annulus braid group in $d$ strings.

\item Moreover, over to each of the two subspaces $\mathbb CP^{d -1}_a$ and 
$\mathbb CP^{d -1}_b$, the map $\Gamma_Q$  inherits a similar stratified structure with respect to the dimensional shift  $d \Rightarrow d - 1$. $\Box$
\end{itemize}
\end{theorem}

\begin{figure}[ht]
\centerline{\includegraphics[height=4in,width=5in]{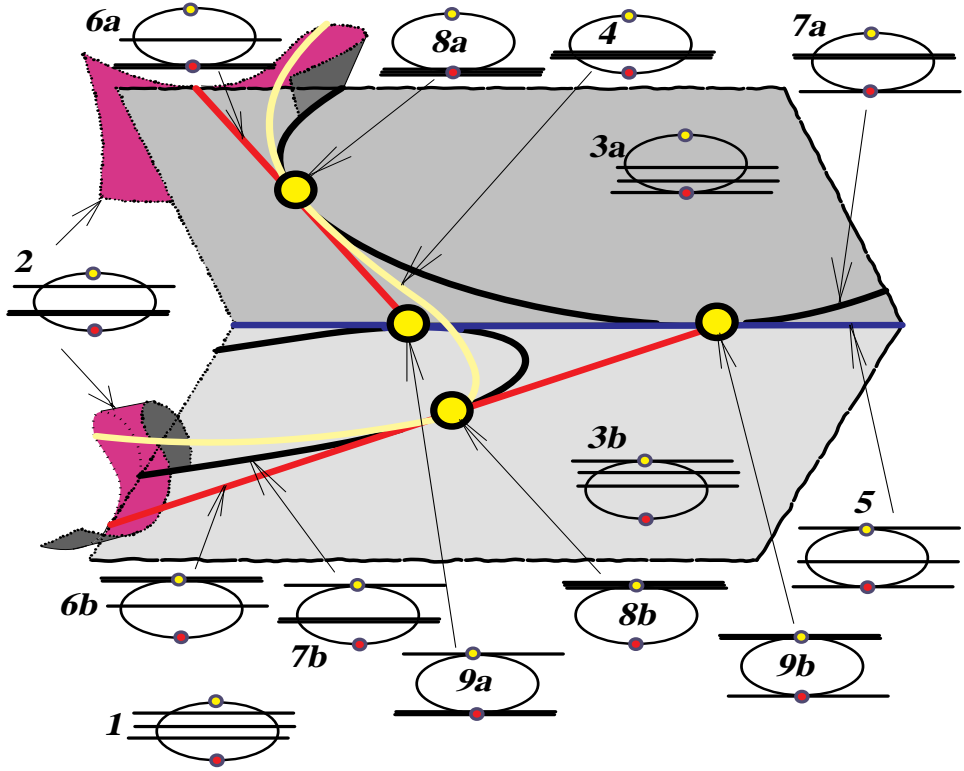}} 
\bigskip
\caption{Patterns labeling the stratification of the map $\Gamma_Q: \mathbb CP^3 \rightarrow CP^3$ and the geometry of its ramification locus.}
\end{figure}

{\bf Example 1. \quad} Let $d = 2$. Then $\Gamma_Q: \mathbb CP^2 \rightarrow \mathbb CP^2$ is of degree $4$. The discriminat parabola $\mathcal D_2$ is given by $\{x^2 - 4yz = 0\}$.  The ramification locus $\mathcal D(\Gamma_Q)$  consists of that parabola together with two tangent lines $\mathbb CP^1_a = \{Ax + y + A^2z = 0\}$ and $\mathbb CP^1_b = \{Bx + y + B^2z = 0\}$ which share  a point 
$E = [-(A +B): AB : 1]$. They are tangent to the parabola at the points $F = [-2A : A^2 : 1]$ and 
$G = [-2B : B^2 : 1]$. Here the parameters $A$ and $B$ depend on the quadratic form $Q(x, y, z)$. The 
$\Gamma_Q$-fibers over $E, F, G$ are singletons; the  cardinality of the fiber over $\mathcal D_2 \setminus (F \cup G)$ is $3$; the cardinality of the fiber over $\mathbb CP^1_a \setminus (F \cup E)$ and over $\mathbb CP^1_b \setminus (G \cup E)$ is  $2$. 
It is still a bit mysterious how all this data manage to produce  $\mathbb CP^2$ as a covering space!  Anyway, the compliment $\mathbb CP^2  \setminus \mathcal D(\Gamma_Q)$ is a $K(\mathsf B_2^{ann}, 1)$-space, where $\mathsf B_2^{ann}$ is the braid group with two strings in a  "fat" annulus. 
\bigskip

Now consider the case $d = 3$ depicted in Figure 3. That figure exhibits a two-fold symmetry that exchanges  points $a$ and $b$ where the pencil of parallel lines is tangent to the curve 
$\mathcal Q$. Pattern 1 in Figure 3 corresponds to the generic stratum $\mathbb CP^3$, while  pattern 2  defines the discriminant surface $\mathcal D_3 \subset \mathbb CP^3$ of degree $4$. It is homeomorphic to $\mathbb CP^1 \times \mathbb CP^1$. Patterns 3a and 3b each  corresponds  to two planes $\mathbb CP^2_a$ and $\mathbb CP^2_b$. Each of the two planes is tangent to the discriminant surface $\mathcal D_3$ along  lines $l_a$ and $l_b$ labeled by patterns 6a and 6b, respectively. These planes intersect along another line $l$ labeled by pattern 5. Pattern 4 encodes the singular locus $\mathcal C$ of the discriminant surface $\mathcal D_3$. It is a rational curve of degree 3 that is homeomorphic to $S^2$. In fact, $\mathcal D_3$ is the ruled surface spanned by lines tangent to $\mathcal C$ ([K]). The surface $\mathcal D_3$ intersects with the plane $\mathbb CP^2_a$ along the union of the line $l_a$ with a quadratic curve $\mathcal D_{2, a}$ labeled by pattern 7a. 
The lines $l_a$ and $l$ are both tangent to the parabola $\mathcal D_{2, a}$. This configuration is already familiar from our description of the case $d = 2$. Similarly,  pattern 7b labels  parabola 
$\mathcal D_{2, b}$. Curves $\mathcal C$, $\mathcal D_{2, a}$, and $l_a$ are all tangent at a point $A_a$ labeled by pattern 8a. Curves $l$ and $\mathcal D_{2, a}$ are tangent at a point $B_a$ labeled by pattern 9b. That point also lies on the line $l_b$. The labeling of points $A_b$ and $B_a$ is done by patterns 8b and 9a, respectively.  

In accordance with Theorem 10,  the ramification locus $\mathcal D(\Gamma_Q)$ for the 8-to-1 map $\Gamma_Q$ coinsides with $\mathcal D_3 \cup \mathbb CP^2_a \cup \mathbb CP^2_b$ and its complement is a $K(B_3^{ann}, 1)$-space. $\Box$

Now, we turn to deconstructions of both homogeneous non-homogeneous polynomials on complex quadratic surfaces $\{Q(x, y, z) = const\}$. First, consider a general homogeneous polynomial $P$ of degree $d$. We can apply  decomposition (6)  to the homogeneous term $R$ of degree $d - 2$ in (6). This will produce a new collection 
of vectors---a new multipole---$\{w_{1, \nu}\}$ associated with an appropriate generalized  parcelling of $\mu : Z(R, Q) \rightarrow \mathbb Z_+$.  
This process of producing lower order multipoles 
$\{w_{s, \nu} \in V(1)\}_{s, \nu}$, where $0 \leq s \leq \lceil d/2\rceil$ and  $0 < \nu \leq d - 2s$, can be repeated again and again until all the  degrees are ``used up". 

We notice that the leading multipole (of highest degree) is determined by this algorithm in a more ``direct way" than the lower degree multipoles. Also note that the choice of a generalized parcelling $\tau$ for 
$\mu: Z(P, Q) \rightarrow \mathbb Z_+$ affects the choice of a generalized parcelling for 
$\mu: Z(R, Q) \rightarrow \mathbb Z_+$, where the polynomial $R = (P - \prod_\nu L_\nu)/Q^s$. Here the 
product $\prod_\nu L_\nu$ is determined by the $\tau$, and  $s$ is the maximal power of $Q$ for which the division in the ring of polynomials is possible.  

{\bf Example 2. \quad} Let  $d = 5$. Then any homogeneous polynomial $P$ of degree $5$ has a representation of the form
\begin{eqnarray}
P  = L_{01}L_{02} L_{03} L_{04}L_{05} + Q \cdot L_{11}L_{12} L_{13}  + Q^2\cdot L_{21} 
\end{eqnarray}
where all the $L_{i j}$'s are linear forms (some of which might be zeros).  The number of such representations does not exceed 
$(9!!) \times (5!!)$. $\Box$

Any non-homogeneous polynomial $P(x, y, z)$ of degree $d$ can be written in the form $P^{(0)} + P^{(1)}$ where the degrees of monomials comprising $P^{(n)}$ are congruent to 
$n$ modulo $2$. Note that the decomposition $P = P^{(0)} + P^{(1)}$ intrinsically makes sense 
on every quadratic surface 
$Q(x, y, z) = \lambda, (\lambda \in \mathbb C)$:  if $P$ belongs to the principle  ideal 
$\langle Q - \lambda \rangle$, so does each term $P^{(n)}, \, n = 0, 1$. Thus, if 
$P \equiv \tilde P\; {\rm mod} \langle Q - \lambda \rangle$, then  we have 
$P^{(n)} \equiv \tilde P^{(n)} \; {\rm mod} \langle Q - \lambda \rangle$. 

On the surface $\{Q(x, y, z) = 1\}$, any component $P^{(n)}$ of the polynomial 
$P^{(0)} + P^{(1)}$ can be homogenized by multiplying its terms of the same degree by 
an appropriate power of $Q$. We denote by $P^{(n)}_Q$ the appropriate homogeneous polynomial.
Generically, $deg(P^{(n)}_Q) = d - n$.

{\bf Example 3. \quad} Let  $d = 5$. Any  polynomial $P$ of degree $5$ has a representation of the form
\begin{eqnarray}
P & = & L_{01}L_{02} L_{03} L_{04}L_{05} + Q \cdot L_{11}L_{12} L_{13}  + Q^2\cdot L_{21} \nonumber\\
& + & M_{01}L_{02} M_{03} M_{04} + Q \cdot M_{11}M_{12}   +  Q^2\cdot  \lambda
\end{eqnarray}
where all the $L_{i j}$'s and $M_{i j}$'s are linear forms (some of which might be zeros), and $\lambda$ is a number. The number of such representations does not exceed 
$(9!!) \times (5!!) \times (7!!) \times (3!!) = 9 \times 7^2 \times 5^3 \times 3^4$. $\Box$

Now one can apply recursively Lemmas 1,   3 and Corollary 5 to each homogeneous polynomial $P^{(n)}_Q$, \, $n = 0, 1$.  Letting $Q = 1$ proves formula $(3)$ from Theorem 2. Let us restate and generalize this theorem in terms of the multipoles:
\begin{theorem} 
\begin{itemize}
\item Let $Q(x, y, z)$ be an  irreducible quadratic form and let $P(x, y, z)$ be any complex polynomial of degree $d$. Its restriction $P|_\mathcal S$ to the complex  surface 
$\mathcal S = \{Q(x, y, z) = 1\}$ admits a representation of the form
\begin{eqnarray} 
P(x, y, z) =  \lambda_0 + \sum_{k = 1}^{d} \; 
\prod_{l = 1}^{k} L_{k , l}(x, y, z),
\end{eqnarray}
where   the linear forms $\{L_{k , l}\}$  are chosen so that each non-zero product 
$\prod_{l = 1}^k L_{k , l}(x, y, z)$ is determined, via the map $\Phi_Q$ 
(see  $(16)$),  by an appropriate multipole  $w_k$ from the variety $\mathcal M(k)$. 

\item The representation $(34)$ is unique, up to a finite ambiguity and up to reordering and rescaling of multipliers  in each  product $\prod_{l = 1}^k L_{k , l}(x, y, z)$. In other words, the set of $P$-representing multipoles $\{w_k \in \overline{\mathcal M}(k)\}_{1 \leq k \leq d}$  is finite. Its cardinality does not exceed 
$\prod_{k = 1}^{d} [(2k - 1)!!]$. \qquad $\Box$
\end{itemize}
\end{theorem}

We would like to end this section by establishing a few facts about  alternative decompositions of polynomials that are based on formula $(9)$. In achieving this goal we are guided by Theorem 22.2 from [Sh].

We claimed that the direct sum in $(7)$ is \emph{orthogonal} with respect to the inner product  $\langle f, g \rangle = \int_{S^2} f\cdot g \; dm$, where $dm$ is the standard rotationally symmetric measure  on the  unit sphere $S^2$. Let us clarify this claim. Applying $(7)$ recursively, we get that for any real homogeneous polynomial $P$ of degree $d$ can be written as $\sum_{d - 2k \geq 0} \; Q^k \cdot P_k^H$, where $P_k^H$ is a real homogeneous harmonic polynomial of degree $d - 2k$. Letting $Q = 1$, we see that for any homogeneous $P$, there is a harmonic polynomial  $P^H = \sum_{d - 2k \geq 0} \; P_k^H$ such that $P|_{S^2} =  P^H|_{S^2}$. Now, any two harmonic homogeneous polynomials $P_k^H$ and $P_l^H$  of different degrees are eigenfunctions with different eigenvalues 
$-(d - 2k)(d - 2k + 1)$ and $-(d - 2l)(d - 2l + 1)$ for the Laplace operator $\Delta_{S^2}$ on the sphere. Therefore, they are orthogonal, i.e. $\int_{S^2} P_k^H \cdot  P_l^H \; dm = 0$. Since any homogeneous polynomial $R$ of degree $d - 2$ can be represented as  $\sum_{d - 2k \geq 0;\;  k > 0} \; Q^k \cdot P_k^H$, we get  the claimed orthogonality $\int_{S^2} (P_0^H)(Q \cdot R) \; dm = 0$ of the direct sum in $(7)$. 

This argument extends to complex harmonic polynomials as follows. We already remarked that the real and imaginary parts of a complex harmonic polynomial $P$ are real harmonic polynomials of the same degree.  If the polynomial is homogeneous, so are its real and imaginary parts. Therefore, for any complex harmonic and homogeneous polynomial $P$ of degree $d$ and any homogeneous polynomial $R$ of degree $d - 2$, we get  $\int_{S^2} P [Q \cdot \bar {R}] \; dm = 
\int_{S^2} (\mathsf{Re} P + \mathbf i\; \mathsf{Im} P) (\mathsf{Re} R -  \mathbf i\; \mathsf{Im} R) \; dm =  
\int_{S^2} (\mathsf{Re} P \cdot \mathsf{Re} R) \; dm +  \int_{S^2} (\mathsf{Im} P \cdot \mathsf{Im} R)\; dm
 -  \mathbf i \int_{S^2} (\mathsf{Re} P \cdot \mathsf{Im} R)\; dm +  \mathbf i \int_{S^2} (\mathsf{Im} P \cdot \mathsf{Re} R)\; dm$. Each of the four integrals must vanish because each integrant is a product of a real homogeneous  and harmonic polynomial of degree $d$ by a real polynomial of a lower degree. 
 
In the spherical coordinates $\theta, \phi$, the inner product is given by an integral 
 $$\int_{S^2} f\cdot \bar g \; dm = \int_{0 \leq \phi \leq \pi;\; 0 \leq \theta \leq 2\pi} f(\Lambda(\theta, \phi))\cdot \bar g(\Lambda(\theta, \phi)) \, |\sin \phi| \; d\theta\, d\phi$$
where $\Lambda(\theta, \phi) = (\cos \theta\, \sin \phi,\; \sin \theta\, \sin \phi,\; \cos \phi)$.

Now we are going to transfer this hermitian inner product from the sphere to its image $\Upsilon_Q$ in a given complex quadratic surface $\mathcal S_Q := \{Q(x, y, z) = 1\}$. As before, let $A$ be a complex invertible matrix that reduces $Q$ to the sum of squares: $x'^2 + y'^2 + z'^2 = Q((x, y, z)\cdot A)$, where $(x', y', z') = (x, y, z)\cdot A$. Let $\Upsilon_Q = \{Q((x, y, z)\cdot A) = 1\}$, with  $(x, y, z)\cdot A$ being a \emph{real} vector.  The ellipsoid  $\Upsilon_Q$ is the image of the unit sphere under the complex linear transformation $A^{-1}$. It is a totally real 2-dimensional algebraic surface in the complex surface $\mathcal S_Q$. 

The real-valued  measure $dm_Q$ on $\Upsilon_Q$ is the pull-back under the map $A$ of the standard measure on the unit sphere. It is invariant under the action of the compact subgroup 
$A\cdot O(3; \mathbb R)\cdot A^{-1} \subset GL(3; \mathbb C)$, where $O(3; \mathbb R)$ denotes the orthogonal group.

For any pair of (homogeneous) complex  functions $f, g$ on $\mathbb C^3$, we get
\begin{eqnarray} 
\int_{\Upsilon_Q} f((x, y, z) A)\cdot \bar g((x, y, z) A) \; dm_Q =  \int_{S^2} f(x', y', z')\cdot 
\bar g(x', y', z') \; dm \nonumber \\ 
= \int_{0 \leq \phi \leq \pi;\; 0 \leq \theta \leq 2\pi} f(\Lambda(\theta, \phi))\cdot \bar g(\Lambda(\theta, \phi)) \, |\sin \phi| \; d\theta\, d\phi
\end{eqnarray}
Applying $(9)$ recursively, gives the following proposition:
\begin{theorem} 
\begin{itemize}
\item The space of complex homogeneous polynomials admits an 
$O_Q(3; \mathbb C)$-invariant decomposition
\begin{eqnarray}
V(d) = \oplus_{d - 2k \, \geq \,  0} \quad Q^k \cdot  Har_Q(d- 2k)
\end{eqnarray}
The summands in  $(36)$ are orthogonal with respect to the hermitian inner product defined by the formula $(35)$. In particular, the homogeneous polynomials $P \in Har_Q(d)$---the solutions of the equation $\Delta_Q(P) = 0$---are  characterized  by the property:   for each $T \in V_Q(d)$, $$\int_{\Upsilon_Q} P\cdot \bar T  \; dm_Q = 0.$$
\item Any polynomial function $F$ on the surface $\mathcal S = \{Q(x, y, z) = const\}$ can be obtained by restricting to $\mathcal S$  a  polynomial $P \in Ker(\Delta_Q)$.\footnote{not necessarily homogeneous, even when $P$ is homogeneous}
\item For any two polynomials $M$ and $N$, ``the complex Dirichlet problem"
\begin{eqnarray}
\big\{ \Delta_Q(P) = M  , \quad P|_{\mathcal S} = N|_{\mathcal S}\big \}
\end{eqnarray}
has a unique \emph{polynomial} solution $P$.
\item Any homogeneous  polynomial $P$ of degree $d$ admits a Maxwell-type representation
\begin{eqnarray}
P(x, y, z) = \sum_{d - 2k\, \geq \,  0} Q(x, y, z)^{d - k  + \frac{ 1}{2}} \cdot 
\nabla_{\mathbf u_{1, k}}
\nabla_{\mathbf u_{2, k}} \dots 
\nabla_{\mathbf u_{d - 2k, k}}
\Big ( Q(x, y, z)^{-\frac{1}{2}}\Big)
\end{eqnarray}

\end{itemize}
\end{theorem}
{\bf Remark \quad} All the statements of Theorem 12, but the last one (dealing with the generalized Maxwell representation $(38)$), can be easily generalized for polynomials in any number of variables.

{\bf Proof \quad} Decomposition $(36)$ follows from $(9)$, and $(38)$ from  $(36)$ together with $(10)$. The second bullet is obviously implied by  $(36)$ and the linearity of $\Delta_Q$. In order to prove the
third bullet, note that $(9)$ implies that $\Delta_Q: V(k) \rightarrow  V(k -2)$ is onto. Therefore, for any $M$, there exists a polynomial $T$, so that $\Delta_Q (T) = M$. On the other hand by bullet two, there exists a harmonic polynomial $R \in Ker(\Delta_Q)$ such that $R|_{\mathcal S} = (N - T)|_{\mathcal S}$. As a result, $P = T + R$ solves the Dirichlet problem.  In order to show that $P$ is unique, it is sufficient to prove that no polynomial of the form $(Q - 1)S$ belongs to $Ker(\Delta_Q)$. Since $\Delta_Q$ is homogeneous of degree $-2$, $\Delta_Q(Q\cdot S) \neq \Delta_Q(S)$, unless $S = 0$. Indeed, let $F$ be the leading homogeneous portion of $S$ of the degree $deg(S)$. Then $\Delta_Q(Q\cdot S) = \Delta_Q(S)$ implies $\Delta_Q(Q\cdot F) = 0$ which is, as we have shown before, impossible, unless $F = 0$. $\Box$


\section{Deconstructing reality}

We have seen that different parcellings of the intersection $Z(P, Q)$ led to different deconstructions of polynomial functions on a quadratic surface $\{Q = 1\}$. The next idea is to invoke symmetry in order to get a unique equivariant parcelling and thus, a \emph{unique} deconstruction of symmetric functions on a quadratic surface supporting an appropriate group action.

We deal with a special, but important case of $\mathbb Z_2$-action. The model example is provided  by the complex conjugation $\tau : (x, y, z) \rightarrow (\overline x, \overline y, \overline z)$. If $Q(x, y, z)$ has real coefficients, both the curve $\mathcal Q \subset \mathbb CP^2_\ast$ and the surface 
$\mathcal S = \{ Q = 1\} \subset \mathbb C^3_\ast$ are invariant under the $\tau$. When $Q$ is a definite real form, the $\mathbb Z_2$-action by conjugation on  $\mathcal Q$ is \emph{free} since the fixed point set $\mathcal Q^{\mathbb Z_2} = \mathcal Q \cap \mathbb RP^2_\ast = \emptyset$. 
For a homogeneous polynomial $P$ with real coefficients the intersection 
$Z(P,  Q) \subset  \mathbb CP^2_\ast$ and the multiplicity function $\mu: Z(P, Q) \rightarrow \mathbb Z_+$ are $\tau$-invariant as well. We notice that a $\tau$-invariant generalized parcelling of such a  $\mu$ will produce a set of $\tau$-invariant  lines $\mathcal L_\nu$. Thus, we can assume that all the linear forms $L_\nu$ in $(6)$ are with real coefficients.
Indeed, given a  representation $P = \lambda\cdot \prod_\nu L_\nu + Q\cdot R$, where $P, Q$ and $L_\nu$'s are real and $\lambda \in \mathbb C$, we get  $P = \overline{\lambda}\cdot \prod_\nu L_\nu + Q\cdot \overline{R}$. Hence, 
$2P = (\lambda + \overline{\lambda})\cdot \prod_\nu L_\nu + Q\cdot (R + \overline{R})$ which implies the validity of $(6)$ over the reals. 

With this observation in mind, we introduce real versions of the multipole spaces (11), (12).  Let
\begin{eqnarray}
 \mathcal M^\mathbb R(k) =  
 \{[\mathbb R^{3 \; \circ}]^k / \Sigma_k^\mathbb R
\end{eqnarray} 
where the group  $\Sigma_k^\mathbb R$ is defined similar to its complex version. The only difference lies in the definition of $H_k^\mathbb R$: it is now a subgroup of 
$(\mathbb R^\ast)^k$, not of $(\mathbb C^\ast)^k$.

The space of real multipoles  $\mathcal M^\mathbb R(k)$ is a space of a principle $\mathbb R^\ast$-fibration over the real variety $\mathcal B^\mathbb R(k) := Sym^k(\mathbb RP^2)$.  We can form an associated line bundle $\eta^\mathbb R(k) = \{\mathcal M^\mathbb R(k) \times _{\mathbb R^\ast} \mathbb R \rightarrow \mathcal B^\mathbb R(k)\}$.
As in the complex case, put $E\eta^\mathbb R(k) = \mathcal M^\mathbb R(k) \times _{\mathbb R^\ast} \mathbb R$ and form the quotient space 
\begin{eqnarray}
\overline{\mathcal M}^\mathbb R(k) = 
E\eta^\mathbb R(k) / \mathcal B^\mathbb R(k)
\end{eqnarray}
by collapsing the zero section to a point $\bf 0$. As before, $\overline{\mathcal M}^\mathbb R(k)$ is a contractible space. 

Real homogeneous  polynomials of degree $d$ form a totally real subspace $V(d; \, \mathbb R)$ in the Hermitian space $V(d)$. The hermitian metric on $V(d)$ defined by $(35)$ generates an eucledian metric on $V(d; \, \mathbb R)$ induced by the inner product 
\begin{eqnarray}
\langle f, g \rangle^\mathbb R \;:=\; \mathsf{Re}\Big\{ \int_{\Upsilon_Q} f(x, y, z)\cdot  \bar g(x, y, z) \, dm_Q \Big\} \nonumber\\ = \;
\mathsf{Re}\Big\{\int_{S^2} f((x, y, z)A^{-1})\cdot \bar g((x, y, z)A^{-1}) \, dm\Big\} 
\end{eqnarray}

Moreover, as in the complex case, the imbedding 
$\beta_Q: V(d - 2; \, \mathbb R)  \rightarrow V(d; \, \mathbb R)$ is an isometry (recall that $Q$  is real).

Similarly to the complex case, we introduce  a variety 
${\mathcal Fact}(d; \, \mathbb R) \subset V(d; \, \mathbb R)$ of completely factorable  real polynomials. 
The orthogonal projection $V(d; \, \mathbb R) \rightarrow V^\perp_Q(d; \, \mathbb R)$ maps ${\mathcal Fact}(d; \, \mathbb R)$ into the vector space $V^\perp_Q(d; \, \mathbb R)$. Due to Theorem 12, $V^\perp_Q(d; \, \mathbb R)$ can be identified with $Ker(\Delta_Q^\mathbb R)$. When the  quadratic form $Q$ is not definite, we can also identify the space  $V^\perp_Q(d; \, \mathbb R)$ with the $d$-graded portion of the polynomial function ring on the real cone $\{Q = 0\}$. 
In any case, we get an algebraic  map
\begin{eqnarray} 
\Phi^\mathbb R_Q : \mathcal M^\mathbb R(d) \rightarrow  {\mathcal Fact}(d; \, \mathbb R) 
 \rightarrow   V^\perp_Q(d; \, \mathbb R)^\circ
\end{eqnarray}
which, by Lemma 1 and the argument above,  is \emph{onto}. 
As in  the complex case, this map extends to  a map 
\begin{eqnarray} 
\tilde \Phi^\mathbb R_Q : E\eta^\mathbb R(d) \longrightarrow  V^\perp_Q(d; \, \mathbb R)
\end{eqnarray}
that takes the zero section $\mathcal B^\mathbb R (d) = Sym^d(\mathbb RP^2)$ of the line bundle $\eta^\mathbb R(d)$ to the origin in $V^\perp_Q(d; \, \mathbb R)$, and each fiber of $\eta^\mathbb R(d)$ isomorphically to a line through the origin. This generates an $\mathbb R^\ast$-equivariant  surjection  
\begin{eqnarray} 
\overline \Phi^\mathbb R_Q : \overline \mathcal M^\mathbb R(d) \longrightarrow  
V^\perp_Q(d; \, \mathbb R)
\end{eqnarray}
with finite fibers, where 
$\overline \mathcal M^\mathbb R(d) := E \eta^\mathbb R (d) / \mathcal B^\mathbb R (d)$.

The number of $\tau$-equivariant  parcellings of $Z(P, Q)$ is harder to determine. It depends only on the restriction of the multiplicity function $\mu$ to the subset $Z(P, Q; \mathbb R) := Z(P, Q) \cap \mathbb RP^2_\ast$ and is equal to the number of generalized parcellings subordinate to such a restriction.  However, if $Z(P, Q; \mathbb R) = \emptyset$, the conjugation acts freely on $\mathcal Q \subset \mathbb CP^2_\ast$ and the $\tau$-equivariant  parcelling is unique.  In such a case, we are getting a more satisfying result:
\begin{theorem}  Let $Q(x, y, z)$ be an irreducible  quadratic form with real coefficients and the signature distinct from $ -3$. Then any real  polynomial $P(x, y, z)$ of degree $d$, being restricted to the real conic 
$\mathcal S^\mathbb R = \{Q(x, y, z) = 1\}$ in $\mathbb R^3_\ast$,  has a  representation 
\begin{eqnarray}
P(x, y, z) = 
\lambda_0 + \sum_{k = 1}^d \lambda_k \Big[\prod_{j = 1}^k  (a_{k, j}x + b_{k, j}y + c_{k, j}z)\Big],
\end{eqnarray}
where each  $(a_{k, j}, \, b_{k, j}, \, c_{k, j})$ is a unit vector and, for $k > 0$,  
$\lambda_k \geq 0$.

When  the surface $\mathcal S^\mathbb R$ is an ellipsoid, the represention $(45)$ is unique, up to reordering and rescaling of the multiplyers in the products. In such a case,  $P$ gives rise to a unique sequence of multiploles 
$\{w_{k} \in \overline{\mathcal M}^{\mathbb R}(k)\}_k$. $\Box$
\end{theorem}
We observe that $\mathcal M^{\mathbb R}(k)$ happens to be a nonsingular space: locally  
$\mathbb RP^2$ and $\mathbb CP^1$ are diffeomorphic, and $Sym^k(\mathbb CP^1) \approx \mathbb CP^k$ is non-singular. An $\mathbb R$-version of the argument that follows Theorem 6 and is centered on formula $(22)$ is valid. Therefore, an $\mathbb R$-version of Lemma 6 holds: the ``combinatorial" and the ``smooth" ramification loci  coincide. 

Note that one always can find  a real homogeneous polynomial $P$ of degree $k$ for which the curves 
$\mathcal P, \mathcal Q \subset \mathbb RP^2_\ast$ have an empty intersection: just take a linear form $L$ such that the line $\mathcal L := \{L = 0\}$ misses the real quadratic curve $\mathcal Q$; then consider any polynomial $P$ sufficiently close to $(L)^k$.  For such a $P$, the $\mathbb Z_2$-equivariant parcelling is unique and, thus, $(\Phi_Q^\mathbb R)^{-1}(P\big|_{Q = 0})$ is a singleton. Evidently, such polynomials $P$ form an open set. Therefore, 
$\Psi_Q^{\mathbb R}: Sym^k(\mathbb RP^2) \rightarrow  \mathbb RP(V^\perp_Q(k; \, \mathbb R))$ is a map of degree one. These remarks lead  to
\begin{theorem}
The multipole space $\mathcal M^{\mathbb R}(k)$, which parametrizes completely factorable real homogeneous polynomials of degree $k$, is a space of an $\mathbb R^\ast$-bundle over the real variety $Sym^k(\mathbb RP^2)$. The $\mathbb R^\ast$-equivariant map 
$\Phi_Q^{\mathbb R}$ induces  a surjective mapping 
$\Psi_Q^{\mathbb R}: Sym^k(\mathbb RP^2) \rightarrow \mathbb RP^{2k} = 
\mathbb RP(V^\perp_Q(k; \, \mathbb R))$ of degree one. Unless $Q$ is a definite form,  the map $\Phi_Q^{\mathbb R}$ is ramified over  the discriminant variety $\mathcal D(\Phi_Q^{\mathbb R})$ of codimension one.  Here $\mathcal D(\Phi_Q^{\mathbb R})$ is comprised of homogeneous polynomials $P$ of degree $k$ on the real surface  $\{Q = 0\}$ for which the surfaces $\{P = 0\}$ and $\{Q = 0\}$ share a line  of multiplicity at least two.\footnote{equivalently, the curves  $\mathcal P$ and $\mathcal Q$ in  $ \mathbb RP^2_\ast$ share a point of multiplicity at least two.}  $\Box$
\end{theorem}
\begin{cor} Let $Q$ be a  positive definite form. Then
\begin{itemize} 
\item the multipole space
 $\mathcal M^\mathbb R(d) =  \Big \{[\mathbb R^{3\; \circ}]^d \Big\}\Big/\Sigma_d^\mathbb R$ 
is diffeomorphic to the space $\mathbb R^{2d + 1\; \circ} \approx S^{2d} \times \mathbb R$.
\item the multipole space $\overline{\mathcal M}^\mathbb R(d)$ in $(40)$ is homeomorphic  to the space $\mathbb R^{2d + 1}$.
\item the real varieties $Sym^d(\mathbb RP^2)$ and $\mathbb RP^{2d}$ are diffeomorphic.
\end{itemize}
\end{cor} 

{\bf Proof \quad} In this case, the ramification locus $\mathcal D(\Phi_Q^{\mathbb R}) = \emptyset$. Since under the corollay's hypotheses the $\tau$-equivariant parcelling is unique, every real homogeneous $P$,  not divisible by $Q$, gives rise to a unique leading multipole $w(P) \in \mathcal M_Q^\mathbb R(d)$---the map $\Phi^\mathbb R$ is $1$-to-$1$. Recalling that $\overline \Phi^\mathbb R$ is onto a vector space of dimension $2d + 1$ and that the smooth ramification locus $\overline\mathcal E = \mathcal D(\overline\Phi_Q^{\mathbb R}) = \emptyset$, completes the proof. $\Box$
\smallskip

The diffeomorphism $Sym^d(\mathbb RP^2) \rightarrow \mathbb RP^{2d}$ serves as yet another  transparent illustration to the Dold-Thom theorem:  not only it reveals $Sym^\infty(\mathbb RP^2)$ as a $K(\mathbb Z_2, 1)$ space, but we actually know how this stabilization of the homotopy groups occurs. In fact, $\pi_1(Sym^d(\mathbb RP^2)) = \mathbb Z_2$,  and for $k >1$,
$\pi_k(Sym^d(\mathbb RP^2)) = \pi_k(S^{2d})$. In particular, $\{\pi_k(Sym^d(\mathbb RP^2))\}$ vanish  for $1 < k < 2d$.
\bigskip

We have seen already the main advantages of the multipole representations for the polynomials on quadratic surfaces: such representations are independent on the  choice of coordinates in $\mathbb C^3_\ast$ or  $\mathbb R^3_\ast$. This is in the sharp contrast with the classical decompositions in terms of the spherical harmonics. In the case of $Q = x^2 + y^2 + z^2$ and over the reals, the independence of the multipoles under the rotations was observed by many. When $Q$ is not positive-definite, similar observations hold.

For a non-degenerated quadratic form  $Q$  with real coefficients,  let  $O_Q(3, \mathbb R)$ denote the group of linear transformations from $GL(3; \mathbb R)$ that preserve the form. When the signature 
$sign(Q) = 2$, then $O_Q(3, \mathbb R)$ contains the Lorenz transformation group (equivalently, the isometry group of a hyperbolic plane) as a subgroup of index two.

The next proposition should be compared with Lemma 2 and Theorem 12.  We  noticed already that, for a \emph{real} quadratic form $Q$, the decomposition (9) is a complexification of a similar decomposition
\begin{eqnarray}
V(d; \mathbb R) = Har_Q(d; \mathbb R) \oplus V_Q(d; \mathbb R)
\end{eqnarray}
over the reals. Here $Har_Q(d; \mathbb R) := Ker(\Delta_Q; \mathbb R)\cap V(d; \mathbb R)$. Therefore, (46) is $O_Q(3, \mathbb R)$-equivariant with respect to the natural action on the space $V(d; \mathbb R)$ and orthogonal with respect to the inner product $(41)$. At the same time, the multipole space $\mathcal M^\mathbb R(d)$ is also equipped with the  $O_Q(3, \mathbb R)$-action induced by the obvious diagonal action on $[V(3; \mathbb R)]^d$. 
Furthermore, the map $\mathcal Fact(d, \mathbb R) \rightarrow Har_Q(d; \mathbb R)$ induced  by the projection $V(d; \mathbb R) \rightarrow Har_Q(d; \mathbb R)$ defined by  (46) is  equivariant.
\begin{cor} Given a real polynomial $P$ on $\mathcal S^\mathbb R$ together with its representation $(45)$ and a transformation $U \in O_Q(3, \mathbb R)$, the transformed polynomial 
$U^\ast(P)(x, y, z) := P((x, y, z)\cdot U)$  on $\mathcal S^\mathbb R$ acquires a representation in the form
\begin{eqnarray}
U^\ast(P)(x, y, z) = 
\lambda_0 + \sum_{k = 1}^d \lambda_k \Big[\prod_{j = 1}^k  (a_{k, j}'x + b_{k, j}'y + c_{k, j}'z)\Big],
\end{eqnarray}
where each new multipole vector $(a_{k, j}',  b_{k, j}', c_{k, j}') = (a_{k, j},  b_{k, j}, c_{k, j})\cdot U^T$. In other words,  the onto map 
$$\Phi_Q^\mathbb R : \prod_{k = 0}^d \mathcal M^{\mathbb R}(k) \rightarrow \prod_{k = 0}^d [V(k, \mathbb R) / V_Q(k, \mathbb R)] \; \approx \oplus_{k = 0}^d \; Har_Q(d; \mathbb R)$$ is 
$O_Q(3, \mathbb R)$-equivariant.   $\Box$
\end{cor}
Combining decomposition $(46)$ with Theorem 12  we get its real analog:
\begin{theorem} Let $Q$ be a  real non-degenerated quadratic  form.
\begin{itemize}
\item The space of real homogeneous polynomials admits an 
$O_Q(3; \mathbb R)$-invariant decomposition
\begin{eqnarray}
V(d ; \mathbb R) = \oplus_{d - 2k \, \geq \,  0} \quad Q^k \cdot  Har_Q( d - 2k; \mathbb R)
\end{eqnarray}
The summands in $(48)$ are orthogonal with respect to the  inner product defined by the formula $(41)$. 
\item Any polynomial function $F$ on the surface $\mathcal S^\mathbb R = \{Q(x, y, z) = const\}$ can be obtained by restricting to $\mathcal S^\mathbb R $  a  polynomial $P \in Har_Q(  \mathbb R)$.
\item For any two polynomials $M$ and $N$, the Dirichlet problem
\begin{eqnarray}
\big\{ \Delta_Q(P) = M  , \quad P|_{\mathcal S^\mathbb R} = N|_{\mathcal S^\mathbb R}\big \}
\end{eqnarray}
has a unique real \emph{polynomial} solution $P$.
\item Any real  homogeneous polynomial $P$ of degree $d$ admits the Maxwell-type representation
\begin{eqnarray}
P(x, y, z) = \mathsf{Re}\Big\{\sum_{d - 2k\, \geq \,  0} Q(x, y, z)^{d - k  + \frac{ 1}{2}} \cdot 
\nabla_{\mathbf u_{1, k}}
\nabla_{\mathbf u_{2, k}} \dots 
\nabla_{\mathbf u_{d - 2k, k}}
\Big ( Q(x, y, z)^{-\frac{1}{2}}\Big)\Big\}
\end{eqnarray}
where $\{\mathbf u_{j, k}\}$ are complex 3-vectors. 
These vectors are real and representation $(50)$ is unique, provided that $Q$ is positive-definite. $\Box$
\end{itemize} 
\end{theorem}


\section{Why one does rarely  see multipoles in non-quadratic skies?}

Let $Q(x, y, z)$ be an irreducible  form of degree $l$ over $\mathbb C$. Then $\mathcal S := \{Q(x, y, z) = 1\}$ is the  surface in $\mathbb C^3_\ast$ and $\mathcal Q := \{Q(x, y, z) =  0\}$ is an irreducible curve in $\mathbb CP^2_\ast$ of degree $l$.

As before, we denote by $V_Q(d)$ the set of homogeneous degree $d$ complex polynomials that are divisible by $Q$. Let $V_Q^\perp(d) \approx V(d)/V_Q(d)$ be an orthogonal complement to $V_Q(d)$ in $V(d)$. 

For any sequence  of non-negative integers $\{d_i\}$ so that $\sum_{1 \leq i \leq s} d_i  = d$, consider a map 
\begin{eqnarray}
\eta:   \prod_{i = 1}^s V_Q^\perp(d_i) \rightarrow V_Q^\perp(d)
\end{eqnarray}
which is defined by taking the product of homogeneous polynomials $P_i \in  V_Q^\perp(d_i)$ restricted to  the surface $\{Q(x, y, z) = 0\}$. 

As before, the subgroup $H_s \subset (\mathbb C^\ast)^s$ of rank 
$s - 1$ acts freely on $ \prod_{i = 1}^s [V_Q^\perp(d_i)^\circ]$ by scalar multiplication. By the definition of $H_s$, the map $\eta$ is constant on the orbits of this action. We view the partition 
$ \{d = \sum_{1 \leq i \leq s} d_i\}$ as a non-increasing function $\omega: i \rightarrow d_i$ on the index set $\{1, 2, \; \dots \;,  s\}$.  Denote by $S_\omega$ the subgroup of 
the permutation group $S_s$ that preserves  $\omega$. Let $\Sigma_\omega$ be an extension of $S_\omega$ by $H_s$ that is generated by the obvious actions of $S_\omega$ and  $H_s$ on $(\mathbb C^\ast)^s \approx \prod_{i = 1}^s \mathbb C^\ast_i$.

Evidently, $\eta$ is an $\Sigma_\omega$-equivariant map. Thus, it gives rise to a well-defined 
map\footnote{To simplify our notations, we do not indicate (as before) the dependency of the map on $Q$.} 
\begin{eqnarray}
\Phi_\omega:   \prod_{i = 1}^s V_Q^\perp(d_i)^\circ \Big / \Sigma_\omega \rightarrow V_Q^\perp(d)^\circ
\end{eqnarray}
Because $Q$ is irreducible, the map $\Phi_\omega$ has  $V_Q^\perp(d)^\circ$, and not just $V_Q^\perp(d)$, as its target. 

As before, we introduce the mulipole space   
\begin{eqnarray}
\mathcal M_Q(\omega) := \prod_{i = 1}^s V_Q^\perp(d_i)^\circ \Big / \Sigma_\omega
\end{eqnarray}
which is  a space of  a principle $\mathbb C^\ast$-fibration over the orbifold 
\begin{eqnarray}
\mathcal B_Q(\omega) : = \prod_{i = 1}^s \mathbb CP(V_Q^\perp(d_i)) \Big /  S_\omega
\end{eqnarray}

As in the case of quadratic forms $Q$, one has a map $\Theta :  \mathcal M_Q(\omega) \rightarrow \mathcal Fact(\omega)$, where $\mathcal Fact(\omega) \subset V(d)^\circ$ is the variety of homogeneous degree $d$ polynomials in $x, y$, and $z$ that admit a factorization as a product of polynomials of the degrees $\{d_i\}_{1 \leq i \leq s}$ prescribed by the partition $\omega$. The map $\Theta$ takes each multipole to the product of the corresponding polynomial factors. Unlike the case of a quadratic $Q, \;$ 
$\Theta$ may not be a 1-to-1 map, although, its generic fiber is a singleton. This conclusion follows from the unique factorization property for the ring 
$\mathbb C[x, y, z]$:  just consider elements of $\mathcal Fact(\omega)$ that are products of \emph{irreducible} factors of the degrees prescribed by $\omega$. The same uniqueness of factorization implies that each fiber of $\Theta$ is finite: there are only finitely many ways of organizing irreducible factors, in which a polynomial $P \in \mathcal Fact(\omega)$ decomposes, into blocks of degrees $\{d_i\}$. 

The map $\Theta$ is not surjective either. However, $\Phi_\omega$ takes the multipole space \emph{onto} the space  ${\mathcal Fact}_Q(\omega)$ of degree $d$ homogeneous polynomials on the surface $\{Q = 0\}$ that admit factorizations subordinate to $\omega$:
\begin{eqnarray} 
\Phi_\omega : \mathcal M_Q(\omega) \stackrel{\Theta}{\rightarrow}  {\mathcal Fact}(\omega)  \stackrel{\Pi}{\rightarrow}  {\mathcal Fact}_Q(\omega) \subset V_Q^\perp(d)^\circ
 \end{eqnarray}
 \begin{definition} Let $d$ be a natural number and $\omega = \{d = \sum_{i = 1}^s d_i\}$ its partition. Let  $Z$ be a finite set equipped with a multiplicity function 
$\mu: Z \rightarrow \mathbb N$ whose $l_1$-norm $\|\mu\|_1$ is 
$ld$.  A \emph{generalized} $\omega$-\emph{parcelling} of $(Z, \mu)$ is a collection of functions 
$\mu_i : Z \rightarrow \mathbb N$, such that 
\begin{itemize}
\item $\sum_i \mu_i = \mu$
\item $\|\mu_i\|_1 = ld_i$
\end{itemize}
When $Z$ is comprised of $ld$ points and each $\mu_i$ takes only two values $0, 1$, the generalized parcelling is called just an $\omega$-\emph{parcelling}.
\end{definition}
Any polynomial $P(x, y, z)$  of degree $d$ that is not divisible by $Q$ determines a multiplicity function  $\mu : Z(P, Q) \rightarrow \mathbb N$  whose $l_1$-norm is $ld$. 
Here $Z(P, Q) := \mathcal P \cap \mathcal Q \subset \mathbb CP^2_\ast$ is a finite set. If such a polynomial $P$ is a product $\prod_i L_i$, where $deg(L_i) = d_i$, then  the $L_i$'s give rise to a 
unique generalized $\omega$-parcelling $\sum_i \mu_i$. 
\begin{lem} Any generalized $\omega$-parceling $\sum_i \mu_i$ of a given multiplicity function $\mu$ on a finite set $Z \subset \mathcal Q$ corresponds to at most one multipole in the space 
$\mathcal M_Q(\omega)$.
 \end{lem}
{\bf Proof \quad} Assume that, for each index $i$,  there exist a polynomial $L_i$ that realizes $\mu_i$ on the finite intersection set $Z \subset \mathcal Q$. Such polynomial is not divisible by $Q$. Put 
 $Z_i = Z(L_i, Q) := \mathcal L_i \cap \mathcal Q$.  Then, employing the Bezout Theorem,  any other polynomial $M_i$ that  realizes the same multiplicity on the same intersection set $Z_i$ must be of the form $M_i = \lambda_i L_i + Q \cdot R_i$. We notice  that if $L_i \in V_Q^\perp(d_i)$, then $M_i \notin V_Q^\perp(d_i)$, provided $R_i \neq 0$. $\Box$
 \begin{cor} The map  $\Phi_\omega$ has finite fibers over ${\mathcal Fact}_Q(\omega)$.
 \end{cor}
 {\bf Proof \quad} An element of $P \in V_Q^\perp(d)^\circ$ is determined, up to proportionality, by  its multiplicity function $\mu: Z(P, Q) \rightarrow \mathbb N$. Now the claim of the corollary follows from Lemma 8 and the observation that a given multiplicity function admits only finitely many generalized $\omega$-parcellings.  $\Box$

Of course, not any generalized $\omega$-parceling on a given pair  $(Z \subset \mathcal Q,\, \mu)$ is realizable by a product  $\prod_i L_i$ with the properties as above. Crudly, this happens because not any  $ld_i$ points on  $\mathcal Q$ can be placed on a curve $\mathcal C_i$ of degree $d_i$ that does not contain $\mathcal Q$ as its component. A curve of degree $d_i$ can always accommodate 
 $(d_i^2 + 3d_i)/2$ points in $\mathbb CP^2_\ast$. Thus, if an inequality $(d_i^2 + 3d_i)/2 \geq d_il$ is valid, that is, if $d_i  \geq 2l - 3$, the right curve might be found; but it is still unclear how to avoid the very real possibility that $\mathcal C_i \supset \mathcal Q$ when $d_i \geq l$.  In fact, such a possibility is a reality! 
 
Unfortunately, unless $l = deg(Q) = 2$ or $s = 1$, the image of $\Phi_\omega$ is of a smaller dimension than the one of the target space $V_Q^\perp(d)^\circ$. 
As a result, \emph{there is no analog of the Sylvester Theorem in the non-quadratic skies}; when $deg(Q) > 2$, a generic element of $V_Q^\perp(d)^\circ$ is \emph{irreducible}. Let us explain these claims.

Recall that for any $d \geq l$,\, 
$dim(V_Q^\perp(d)) = \frac{1}{2}\{(d^2 + 3d) - [(d - l)^2 + 3(d -l)]\} = \frac{l}{2}(2d - l +3)$. Since the dimension of the group $H_s$ is $s -1$, we get 
$$dim(V_Q^\perp(d)) - dim(\mathcal M_Q(\omega)) = 
\frac{l}{2}(2d - l +3)  - \sum_{i = 1}^s \frac{l}{2}(2d_i - l + 3) + (s - 1) = 
\frac{s - 1}{2}[l^2 - 3l + 2],$$  provided all $d_i \geq l$.   Under these hypotheses, the difference of the two dimensions vanishes only when $l = 1, 2$ or $s = 1$. Hence, 
\begin{lem}
For $l = deg(Q) > 2$ and $\{d_i \geq l\}_{1 \leq i \leq s}$,  the map $\Phi_\omega$ is not onto, i.e.  a generic polynomial from $V_Q^\perp(d)^\circ$ is not a product of  polynomials  of degrees $\geq l$. The codimension of $\Phi_\omega(\mathcal M_Q(\omega))$ in $V_Q^\perp(d)$ is 
$\frac{s - 1}{2}[l^2 - 3l + 2]$.
\end{lem}
For example, on a cubic surface , $dim(V_Q^\perp(d)) - dim(\mathcal M_Q(\omega)) = 
 (s - 1)$, provided  $\{d_i \geq 3\}_{1 \leq i \leq s}$.
 
If we drop the hypotheses $\{d_i \geq l\}$, the computation is a bit more involved: 
$$ dim(V_Q^\perp(d)) - dim(\mathcal M_Q(\omega)) = 
\frac{l}{2}(2d - l +3)  - \sum_{i : \; d_i \geq l } \frac{l}{2}(2d_i - l +3) 
- \sum_{j : \; d_j <  l} \frac{d_j}{2}(d_j + 3) + (s - 1)$$

We conjecture that the RHS of the formula above is always positive, unless $l = 1, 2 $ or $s =1$.


\section{Multipoles and  function approximations  on quadratic surfaces}

We will be concerned with polynomial approximations of holomorphic functions 
 $f : \mathcal S \rightarrow \mathbb C$  on   an irreducible complex quadratic surface 
 $\mathcal S_Q = \{Q(x, y, z) = 1\}$, as well as  with polynomial approximations of continuous functions  
 $f : \mathcal S^\mathbb R_Q \rightarrow \mathbb R$ on its real version. 
 Eventually, we would like to  represent the approximating polynomials in terms of  their multipoles. 
 
 As one uses polynomials of higher and higher degrees to improve approximations, the issue is  stability of the multipole representations.  In general, such stability is absent for several reasons: 1) the intrinsic ambiguities of the multipole representations for complex plynomials; 2) the impossibility of converting an analytic function on a surface into a ``homogeneous" analytic function in the ambient space (homogenizing polynomials worked well). Even abandoning multipole representations in favor of linear methods of harmonic analysis, does not eliminate the stability issue instantly: in general, 
the coefficients of approximating polynomials fail to stabilize. However, if the approximating polynomials are linear combinations of mutually \emph{orthogonal} and normalized  polynomials (analogous to the  Legendre polynomials), the coefficients of the combinations will stabilize. By introducing  appropriate notions of orthogonality for polynomials on a quadratic surface, we aim to establish similar facts for polynomial approximations there (the spherical harmonics reflect  a particular case of such orthogonality). Then we can combine harmonic analysis with non-linear methods of multipole representation for polynomials. We  call such an approach \emph{synthetic}.  \smallskip

Let $C(K)$ denote the algebra of all continuous $\mathbb C$-functions on a Hausdorff compact space $K$. Recall that a \emph{uniform algebra}  is a closed (in the $sup$-norm) subalgebra $\mathcal A \subset C(K)$  that separates points of $K$. Such an algebra is called \emph{antisymmetric}, if any real-valued function from $\mathcal A$ is constant. A subset $Y \subset K$ is called an \emph{antisymmetry set} for $\mathcal A$ if any function from  $\mathcal A$, which is real-valued on $Y$, is a constant. The Bishop Theorem about antisymmetric subdivisions (cf.  [G], Theorem 13.1) claims that the maximal sets of antisymmetry $\{E_\alpha\}$ are closed and disjoint,  and their union is $K$. Moreover, if $f \in C(K)$ and $f \big|_{E_\alpha} \in \mathcal A  \big|_{E_\alpha}$, then $f \in \mathcal A$. In particular, if each $E_\alpha$ is a singleton, then $\mathcal A = C(X)$. 

For a  space $X \subset \mathbb C^n$, let us denote by $\bar \mathcal P(X)$ the closure in the 
$sup$-norm on compacts in $X$ of the algebra $\mathcal P(X)$ generated by all complex polynomial functions. 
Note that if any two points in $X$ can be separated by a \emph{real-valued} polynomial, then Bishop's Theorem implies that $\bar \mathcal P(X) = C(X)$.

For instance, consider  a section $H$ of the complex surface $\mathcal S_Q$ by a totally real subspace $V^3 \subset \mathbb C^3$ --- an image of $\mathbb R^3 \subset \mathbb C^3$ under a complex transformation $A \in GL(3; \mathbb C)$. Since any two points in $\mathbb R^3$,  can be separated by a real-valued polynomial, the same property holds for any two points in $V^3$, and thus, in $H$. 
By the Bishop Theorem, any continuous function $f$ on $H$ admits an approximation in the $sup$-norm on compacts by complex polynomials. In particular, this conclusion is valid for the real surfaces 
$H = \Upsilon_Q \subset (\mathbb R^3)A^{-1}$ and $H = \mathcal S^\mathbb R_Q \subset \mathbb R^3$ which have been employed on many occasions.

For a compact set $K \subset \mathcal S_Q$, denote by $\hat K$ the polynomial hull (closure) of $K$. It consists of all points $v$ in $\mathbb C^3$ with the property:  $|P(v)| \leq sup_{w \in K} |P(w)|$ for 
\emph{any} complex polynomial $P$. Because for any point $v \notin  \mathcal S$ and $w \in \mathcal S$, we have $|(Q - 1)(v)| > |(Q - 1)(w)| = 0$, the polynomial closure $\hat K$ must be contained in $\mathcal S$. 
In fact, the polynomial closure $\hat K$ of a compact $K \subset \mathbb C^3$ must be a polynomially convex compact set.
A theorem by Oka and Weyl (cf.  [G], Theorem 5.1) claims that any complex analytic function, defined in a neighborhood of a compact polynomialy convex set, admits an approximation in the $sup$-norm on $K$ by complex polynomials. Thus,  any analytic function $f(x, y, z)$ defined in a neighborhood of $\hat K \subset \mathcal S_Q$, $K$ being a compact in $\mathcal S_Q$, can be uniformly approximated on $\hat K$ by complex polynomials. 

When dealing with families of functions on non-compact sets, the default topology in the relevant functional spaces is defined by  the uniform convergency on compact subsets. In this topology, the subalgebra $\mathcal O(\mathcal S_Q)$ of holomorphic functions is closed in the algebra of all continuous functions $C(\mathcal S_Q)$ (cf. [GuR], Lemma 11). In particular, if a sequence of polynomials (in $x, y$, and $z$) is converging in the $sup$-norm on every compact in $\mathcal S_Q$, then its limit is a holomorphic function on $\mathcal S_Q$. One can employ any expanding family $\{K_r\}_{1 \leq r \leq \infty}$ 
(i.e., $K_r \subset K_{r + 1}$ and $\cup_r \; K_r = \mathcal S_Q$)  of polynomially convex compacts $K_r \subset \mathcal S_Q$ to build a sequence $\{P_r\}$ of polynomials that will approximate a given holomorphic function $f \in \mathcal O(\mathcal S_Q)$. 

Let $\mathcal F(\mathcal S_Q)$ be a subset of  $\mathcal O(\mathcal S_Q)$ formed by  functions on 
$\mathcal S_Q$ that admit a representation as a series 
\begin{eqnarray}
\sum_{k = 0}^\infty \; \prod_{j = 1}^k L_{kj},
\end{eqnarray}
where each $L_{kj}$ is a linear form in $x, y$, and $z$. The series is required to converge uniformly  on each compact $K \subset \mathcal S_Q$ (as we remarked before, any such uniformly converging series produces a holomorphic function on $\mathcal S_Q$). 

In fact, one can define a similar set $\mathcal F(K) \subset C(K)$ for any closed $K \subset \mathbb C^3$. Due to Theorem 11, any polynomial on $\mathcal S_Q$ is of the form $(34)$ (which is a special case of $(56)$). Therefore, when $\overline \mathcal P(K) = C(K)$ for a compact 
$K \subset \mathcal S_Q$, then $\mathcal F(K)$ is dense in $C(K)$ as well. 

Examining $(56)$, we observe that if this series converges at a point $v = (x, y, z)$ it  must converge absolutely at any other point $\lambda \cdot v$, were the complex number $\lambda$ has modulus less than 1. For any set $Y \subset \mathbb C^3$, we denote by $Y^\bullet$ the set $\{ \lambda v \big | \; v \in Y, \lambda \in \mathbb C,  |\lambda| < 1\}$ and call it the \emph{round hull} of $Y$.

Consider the set $\mathcal S_Q^\bullet$. Because any complex line through the origin that does not belong to the complex cone $\{Q = 0\}$ hits $\mathcal S_Q$ at a pair of antipodal points, $\mathcal S_Q^\bullet$ is an open domain in $\mathbb C^3$, complementary to the  cone (the origin belongs to $\mathcal S_Q^\bullet$), whose boundary contains $\mathcal S_Q$. Therefore, any function 
$f \in \mathcal F(\mathcal S_Q)$ must be a restriction of a function which is \emph{analytic} in $\mathcal S_Q^\bullet$ and continuous in $\mathcal S_Q^\bullet \cup \mathcal S_Q$. I doubt that the converse  statement  is true. Note that a given function $f$ on $\mathcal S_Q$ may have many analytic extensions in $\mathcal S_Q^\bullet$: for example, 1 extends to 1 and to $1/Q$ (the latter  has poles along the complex cone). 

If a  section $H \subset \mathcal S_Q$ by a totally real subspace $V^3$  is an ellipsoid, then its interior in $V^3$ coincides with $\mathcal S_Q^\bullet \cap V^3$. Thus, any real function $f$ on $H$ that admits a representation as in $(56)$ must be real analytic in the interior of the ellipsoid. So, it is represented by its Taylor series at the origin; on the other hand, series $(56)$, uniformly converging in the vicinity of the origin, is a form of a very specialized Taylor series (just count the dimensions of the coefficient spaces of each degree to see how special it is).

Since any function from  
$\mathcal O(\mathcal S_Q)$ admits a  polynomial approximation on compacts, the subset 
$\mathcal F(\mathcal S_Q)$ is dense in in the space of all holomorphic functions. So, 
$\mathcal F(\mathcal S_Q)$ is squeezed between the vector space $\mathcal O(\mathcal S_Q)$ and its dense subspace $\mathcal P(\mathcal S_Q)$. It is not even clear whether $\mathcal F(\mathcal S_Q)$ is a vector space. To understand the structure of  the set $\mathcal F(\mathcal S_Q)$ seems to be an interesting and hard problem. Unfortunately, our progress towards this goal is minimal.

Note that,  if a holomorphic function vanishes on a totally real analytic surface 
$H \subset \mathcal S_Q$, it must vanish in a neighborhood of $H$ in $\mathcal S_Q$.

We summarize these observations in the following proposition that, in particular, generalizes the second claim in Theorem 3.
\begin{theorem}
Any holomorphic function on  complex surface $\mathcal S_Q$ is a limit in the topology of uniform convergence on compacts  of polynomial functions in $\mathbb C^3$; the approximating polynomials each admit a representation as in $(34)$.

A subset $\mathcal F(\mathcal S_Q)$ of  functions that admit a representation as in $(56)$ is dense in the space of all holomorphic functions $\mathcal O(\mathcal S_Q)$ and invariant under multiplications by scalars  and the natural $O_Q(3; \mathbb C)$-action. Each function $f$ from $\mathcal F(\mathcal S_Q)$ admits a canonic holomorphic extension 
$$f (\lambda x, \lambda y, \lambda z) = \sum_{k = 0}^\infty \; \lambda^k \Big[\prod_{j = 1}^k  L_{kj}(x, y, z)\Big]$$
into the round hull $\mathcal S_Q^\bullet$. Here $(x, y, z) \in \mathcal S_Q$ and the module of $\lambda$ is less than one.

Any continuous function $f: H \rightarrow \mathbb C$ on the totally real  surfaces  $H = \Upsilon_Q$ or 
$H = \mathcal S^\mathbb R_Q$ is a limit in the topology of uniform convergence on compacts in $H$ of polynomial functions in $\mathbb C^3$. As a result, the set $\mathcal F(H)$ is dense in $C(H)$. Each function from $\mathcal F(H)$ admits a similar analytic extension into the open portion of real  cone over $H$ that  is bounded by $H$ and contains the origin.

If $f = F\big |_{H}$ for a  function $F$,  holomorphic in a neighborhood of $H$ in $\mathcal S_Q$, then such an  $F$ is unique.  $\Box$
\end{theorem}

Now consider the vector space $\mathcal P(\mathcal S_Q)$ of all polynomial functions restricted to $\mathcal S_Q$ and equipped  with the inner product $\langle f, g \rangle$ defined by $(35)$. If a \emph{homogeneous} polynomial $P \neq 0$, then $\langle P, P \rangle > 0$. An homogeneous polynomial is determined by its restriction to $\mathcal S_Q$.  Evidently, $\langle P, P \rangle = 0$ implies that the restriction of $P$ to $\Upsilon_Q$ is zero. Since $\Upsilon_Q \subset \mathcal S_Q$ is a totally real analytic submanifold, $P|_{\Upsilon_Q} = 0$ implies that $P$ vanishes in the vicinity of $\Upsilon_Q$ in 
$\mathcal S_Q$. By analyticity, $P$ must vanish everywhere in $\mathcal S_Q$, and hence,  $P$ is a zero polynomial. As a result, $\langle P, P \rangle$ gives rise to a norm on the space of homogeneous polynomials $V(d)$. Because any function from $\mathcal P(\mathcal S_Q)$ is a restriction of an homogeneous polynomial, we get a non-degenerated Hermitian inner product on the vector space $\mathcal P(\mathcal S_Q)$.

In view of Theorem 12 and by a similar line of arguments, $\langle P, P \rangle = 0$ implies that $P = 0$ for any complex polynomial $P \in Ker(\Delta_Q)$. Therefore, being restricted to a space of $Q$-harmonic polynomials,  the inner product $\langle f, g \rangle$ in $(35)$ gives rise to an Hermitian structure and an $L_2$-norm $\|P\|_{\Upsilon_Q}$. In particular, each space $Har_Q(k) \approx V_Q^\perp(k)$ inherits this norm, and  $Har_Q(k)$ is orthogonal to  $Har_Q(l)$, provided $l \neq k$. 

Consider the vector space $\prod_{k = 0}^\infty \; Har_Q(k)$ formed by finite sequences of vectors $\{P_k \in Har_Q(k)\}$ and its closure $L_2^{Har} : = \oplus_{k = 0}^\infty \; Har_Q(k)$ formed by  infinite sequences 
$\{P_k \in Har_Q(k)\}$ subject to the condition $\sum_{k = 0}^\infty \|P_k\|^2_{\Upsilon_Q} < \infty$.

Let $L_2(\Upsilon_Q)$ denote the complex Hilbert space of $L_2$-integrable functions on the ellipsoid $\Upsilon_Q$. Every function $f \in L_2(\Upsilon_Q)$ defines a unique system of its ``Fourier components" $\{f_k \in  Har_Q(k)\}$. Each $f_k$ is the unique polynomial from $Har_Q(k)$ that delivers the minimum $min_{P \in Har_Q(k)} \; \|f  - P\|_{\Upsilon_Q}$. By Theorem 19, any 
$f \in C(\Upsilon_Q)$ is a limit in the $sup$-norm on $\Upsilon_Q$ of harmonic polynomials, it must be also the limit of the same  sequence of polynomials viewed as elements of  $L_2(\Upsilon_Q)$. Therefore, if $f \in C(\Upsilon_Q)$ is orthogonal to all the subspaces $Har_Q(k)$, it must be a zero function. Hence, as an element of $L_2(\Upsilon_Q)$, $f = \sum_{k = 0}^\infty \; f_k$, and 
$\|f \|^2_{\Upsilon_Q} = \sum_{k = 0}^\infty \; \|f_k \|^2_{\Upsilon_Q}$. In particular, any $f \in \mathcal O(\mathcal S_Q)$ determines a unique system of its harmonics $\{f_k \in Har_Q(k)\}$. Moreover, it is determined by 
$f \big |_{\Upsilon_Q}$, and thus, by its harmonics $\{f_k\}$. 

Evidently, the sequence of partial sums $\{f_{[d]} := \sum_{k = 0}^d f_k \in Ker(\Delta_Q)\}_d$ converges in $L_2(\Upsilon_Q)$ to $f$.
In terms of homogeneous polynomials, we get an analogous sequence 
$\{f_{\{d\}} := \sum_{k = 0}^d Q^{\lceil(d - k)/2\rceil}f_k \}_d$ converging  in $L_2(\Upsilon_Q)$ to $f$.

Similar arguments can be applied, instead of  the ellipsoid $\Upsilon_Q$, to any surface $H \subset \mathcal S_Q$ that is a section of $\mathcal S_Q$ by a totally real subspace $V^3 \subset \mathbb C^3$. In particular, they are valid for $\mathcal S_Q^\mathbb R$.  First, we pick a measure $dm$ on $H$ such that any polynomial $P|_H \in L_2(H)$ (when $H$ is compact, one can choose any measure of finite volume).  For example, consider the area 2-form $\omega$ on $H$ induced by the imbedding $H \subset \mathbb C^3$ and multiply it by the function $exp[-(x\bar x + y\bar y + z\bar z)]$ to get the right $dm$. Then we define a new inner product by 
$\langle f, g \rangle_H = \int_H \; f \cdot \bar g \; dm$. Since $H$ is totally real,
this will give rise to an Hermitian structure in each space $V(k)$. Notice that the multiplication-by-$Q$ imbedding $V(k) \rightarrow V_Q(k + 2)$ is an isometry. As before, we form the orthogonal compliments $V_Q^\perp(k)$ to the subspaces  $V_Q(k)$ (the only difference is that now  the space $V_Q^\perp(d)$ could be different from the space $Har_Q(d)$ of $Q$-harmonic polynomials). Then we use $\prod_{k = 0}^\infty \; V_Q^\perp(k)$ to construct  a Hilbert space $\oplus_{k = 0}^\infty \; V_Q^\perp(k)$. As before, any continuous function $f$ on $H$  will acquire a unique a representation 
$f = \sum_{k = 0}^\infty \; f_k$, where $\{f_k \in V_Q^\perp(k)\}$, and $\|f \|^2_H = 
\sum_{k = 0}^\infty \; \|f_k \|^2_H$.

Thanks to Theorems 7 and 11, each complex homogeneous polynomial   $f_k \in V_Q^\perp(k)$ can be represented  by some  multipole $w_k^f \in \overline {\mathcal M}(k)$. Similarly, any real homogeneous polynomial $f_k \in V_Q^\perp(k; \mathbb R)$ can be represented  by some  multipole $w_k^f \in \overline {\mathcal M}^\mathbb R(k)$.
When $Q$ is positive definite, these real multipoles are unique. However, in general,  the ambiguity of the multipole representation could cause some trouble. So, we need to choose the representing multipoles with some care. \smallskip

Due to the embedding  $\Theta: \overline{\mathcal M}(k) \rightarrow V(k)$ (see (16)) with the image $\mathcal Fact(k)$, the multipole space $\overline{\mathcal M}(k)$ acquires a metric $\rho$ induced by the $L_2^H$-norm $\|P\|_H$ in $V(k)$. 

The lemma  below helps to estimate the size of the fiber $\overline{\Phi}_Q^{-1}(u)$ over $u \in V_Q^\perp(k)$ in terms of an universal angle $\theta_k = \theta(k, H, dm)$ and the $L_2^H$-norm of $u$:
\begin{lem} Consider the distance function $\rho$ on $\overline{\mathcal M}(k)$ that is induced by the $L_2^H$-norm $\|\sim\|_H$ in $V(k)$ via the embedding $\Theta: \overline{\mathcal M}(k) \rightarrow \mathcal Fact(k) \subset V(k)$. 
Then there exist an angle
$0 < \theta_k \leq \pi/2$ so that, for each $u \in V_Q^\perp(k) \subset V(k)$, the distance from any 
$w \in \overline{\Phi}_Q^{-1}(u)$ to the zero multipole is at most $\|u\| /sin(\theta_k)$, and the diameter  of the fiber $\overline{\Phi}_Q^{-1}(u)$ is at most $2\|u\| /sin(\theta_k)$. 
\end{lem}

{\bf Proof \quad} Let $S(k)$ denote a unit sphere (with respect to the $\|\sim\|_H$-norm) in $V(k)$ and centered at the origin. Because $Q$ is irreducible, $\mathcal Fact(k) \cap V_Q(k) =  \emptyset$.  Thus, the compact sets $S(k) \cap \mathcal Fact(k)$ and $S_Q(k) = S(k) \cap V_Q(k)$ are disjoint.  Therefore, there is a number $0 < \theta \leq \pi/2$ so that the angle between any two vectors $u \in  \mathcal S(k) \cap \mathcal Fact(k)$ and $v \in S_Q(k)$ is greater than or equal to $\theta$. Note that $S(k) \cap \mathcal Fact(k)$ and $S_Q(k)$ are invariant under the circle action $S^1 \subset  \mathbb C^\ast$. So, 
$\mathcal Fact(k)$ and $V_Q(k)$ also are real cones with their tips at the origin and bases  
$\mathcal S(k) \cap \mathcal Fact(k)$ and  $S_Q(k)$. We conclude that the angle between any two vectors $u \in \mathcal Fact(k)$ and $v \in V_Q(k)$  has the same lower bound $\theta > 0$. Now consider an open real cone $C_Q(k) \subset V(d)$ comprised of  vectors that form an angle $\phi < \theta$ with the subspace $V_Q(k)$ and a  complementary cone $C_Q^\perp(k) :=  V(d) \setminus C_Q(k) \supset V_Q^\perp(k)$. The argument above shows that $\mathcal Fact(k) \subset C_Q^\perp(k)$. Hence, the distance from any $w \in \overline{\Phi}_Q^{-1}(u)$ to the zero multipole is at most $\|u\| /sin(\theta_k)$.  As a result, the diameter  of the fiber $\overline{\Phi}_Q^{-1}(u)$ is at most $2\|u\| /sin(\theta_k)$. $\Box$
\begin{cor} Consider a continuous function $f \in L_2(H)$ and its orthogonal decomposition  \hfill\break $f = \sum_{k = 0}^\infty f_k$, where $f_k \in V_Q^\perp(k)$. Then, for any choice of the multipoles 
$w_k \in \overline{\Phi}_Q^{-1}(f_k)$, 
\begin{eqnarray}
 \sum _{k = 0}^\infty  sin^2(\theta_k)\cdot \rho(w_k, \mathbf 0)^2 < \infty,
\end{eqnarray} 
 where $\rho(w_k, \mathbf 0) =  \| \Theta(w_k) \|_H$.      \quad $\Box$
\end{cor}

It seems to be  far from trivial to understand the asymptotic behavior of  $\{sin(\theta_k)\}$ as 
$k \rightarrow \infty$. Perhaps, the lack of understanding of this asymptotics  it is the most significant gap in our analysis.

To state the last claim in the next theorem, we need one technical definition that  likely has very little to do with the essence of the  statement. The set  $\mathcal D(\Phi_Q)  \subset V_Q^\perp(k)$ is a complex algebraic variety that is stratified by algebraic sets $\{\mathcal D_{k, \pi}\}$ which are labeled by various partitions $\pi$ of $2k$. This labeling is done by attaching the divisor $\mathcal P \cap \mathcal Q \in Sym^{2k}(\mathcal Q)$, or rather the partition $\pi$ of $2k$ defined by the multiplicity function of $\mathcal P \cap \mathcal Q$, to each homogeneous polynomial $P(x, y, z)$ restricted to the cone $\{Q = 0\}$ and viewed as an element of $V_Q^\perp(k)$. In particular, when  
$\pi = \{2d = 1 + 1 + 1 + \dots + 1\}$ or $\{2d = 2 + 1 + 1 + \dots + 1\}$,  then
$\mathcal D_{k, \pi} = V_Q^\perp(k)$ or $\mathcal D_{k, \pi} = \mathcal D(\Phi_Q)$, respectively. The natural partial order among partitions reflects the inclusions of the corresponding strata. If we delete all the substrata from a given stratum $\mathcal D_{k, \pi}$, we get a ``pure" stratum that we denote $\mathcal D_{k, \pi}^\circ$. The  variety $\mathcal D(\Phi_Q)$  is a Whitney stratified space; as a result, the vicinity of every stratum $D_{k, \pi}^\circ$ has a structure of a bundle whose fiber is a real cone over another stratified space $Lk_{k, \pi}$ ---the link of  $D_{k, \pi}^\circ$. We will make use of this fact together with another important feature of the stratification $\mathcal D_{k, \pi}$: namely, all the strata of $Lk_{k, \pi}$ have  even real codimensions.

We say that a parametric curve $\gamma: [0, 1] \rightarrow V_Q^\perp(k) \approx Har_Q(k)$ is \emph{tame} if it consists of a finite number of arcs, each of which has the following property:   the interior of each arc is contained  in some  stratum $\mathcal D_{k, \pi}^\circ$.
We say that a continuous function family $\{f_t \in L_2(\Upsilon_Q)\}_{0 \leq t \leq 1}$  is \emph{tame}, if for each $k$, the path 
$\{(f_t)_k \in Har_Q(k)\}_{0 \leq t \leq 1}$ is tame. 

 First, consider the functions $f \in \mathcal O(\mathcal S_Q)$  such that, for each $k$,  the polynomial 
$f_k \in Har_Q(k) \approx V_Q^\perp(k)$ does not belong to the ramification locus 
$\mathcal D(\Phi_Q) \subset V_Q^\perp(k)$ of the map $\Phi_Q$ (the rest of the functions form a complex codimension one subset $\mathcal D \subset \mathcal O(\mathcal S_Q)$). Over the compliment to the  variety $\mathcal D(\Phi_Q)$, the map $\Phi_Q$ is a covering map. Thus, for each initial lifting, a deformation $(f_t)_k$ admits a unique lifting to the multipole space, as long as   $(f_t)_k \in V_Q^\perp(k) \setminus \mathcal D(\Phi_Q)$. 

Next, for any tame $t$-family $\{f_t \in  \mathcal O(\mathcal S_Q)\}_{0 \leq t \leq 1}$, consider the tame curve 
$\{(f_t)_k \in  V_Q^\perp(k)\}_{0 \leq t \leq 1}$ and the first arc $\gamma$ in a finite sequence of arcs that form this curve. There are two possibilities: 1) the arc starts at a stratum $\mathcal D_{k, \pi}^\circ$ and is confined to it for a while, 2) the arc starts at a stratum $\mathcal D_{k, \pi}^\circ$ but  moves instantly into an ambient stratum $\mathcal D_{k, \pi'}^\circ$. In the first case, over $\mathcal D_{k, \pi}^\circ$, $\Phi_Q$ is a covering map and there is a unique lifting of $\gamma$ extending each lifting $\tilde\gamma(0)$ of $\gamma(0)$. In the second case, we claim  that, for any lifting $\tilde\gamma(0)$ of $\gamma(0)$, in the vicinity of  $\tilde\gamma(0)$ the  map $\Phi_Q$ is onto. Indeed, it is a  proper holomorphic map, and thus, its image must be an analytic space (see [N], Theorem 2, page 129). Because $\Phi_Q$ is finite, the image of a neighborhood of $\tilde\gamma(0)$ under $\Phi_Q$ must be of the maximal dimension, and hence, must contain a neighborhood of  $\gamma(0)$. As a result,  the set  $\Phi_Q^{-1}(\gamma)$ (it is a finite graph) must be present in any neighborhood of  $\tilde\gamma(0)$; so, we can lift $\gamma$ to  an arc $\tilde\gamma$ that starts at $\tilde\gamma(0)$. An induction by the number of arcs in the curve $(f_t)_k$ proves the existence  of its lifting to the multipole space. Note that an analogous argument fails for $\Phi_Q^\mathbb R$: a finite image of a real analytic set is a real semi-analytic set that can miss the arc $\gamma$. Therefore, in Theorem 22 the lifting property for tame deformations is absent. 

The arguments above prove the following theorem:

\begin{theorem} Let $Q(x, y, z)$ be an irreducible complex  quadratic form, and let
$\mathcal S_Q$ be a complex quadratic surface defined by the equation $\{Q = 1 \}$. Denote by $A$ a complex change of coordinates that reduces the form $Q$ to the sum of squares. Let  $\Upsilon_Q \subset \mathcal S_Q$ be a totally real  ellipsoid   defined by the equations  $\{Q((x, y,z)A) = 1, \; \mathsf{Im}((x, y,z)A) = 0\}$ and equipped with the measure defined by $(35)$. Let  $f$ be an analytic function on $\mathcal S_Q$.  
 
 Then there exists sequence  of multipoles 
 $\{w_k^f \in \overline {\mathcal M}(k)\}_{0 \leq k < +\infty}$   such that:
 \begin{itemize}
 \item the sequence  gives rise, via the maps $\{\overline \Phi_{Q}\}$,  to complex $Q$-harmonic polynomials $P^f_d(x,y,z) = \sum_{k = 0}^d \overline \Phi_{Q}(w_k^f)$, where the mutually orthogonal polynomials $\{\overline \Phi_{Q}(w_k^f) \in Har_Q(k)\}$ are uniquely determined by $f$.
 \item as $d \rightarrow \infty$,   the polynomials $\{P_d^f\}$ converge in the space $L_2(\Upsilon_Q)$ to the function $f\big|_{\Upsilon_Q}$, and therefore, uniquely determine $f \in \mathcal O(\mathcal S_Q)$.
\item  for a given $f$,  there are at most $(2k - 1)!!$ choices for  each multipole $w_k^f $. 
\item the multipoles  $\{w_k^f \}$ satisfy  property $(57)$ from Corollary $20$.
\item for any  tame deformation $\{f_t \in \mathcal O(\mathcal S_Q)\}_{ 0 \leq t \leq 1}$ of the function 
$f = f_0$, there  exists a continuous deformation  $\{w_k^{f_t} \in \overline {\mathcal M}(k)\}$ of the $f_t$-representing  multipoles, such that $\{w_k^{f_0} = w_k^f\}$.  For functions $f$ and their continuous deformations $f_t$ outside a subspace $\mathcal D \subset \mathcal O(\mathcal S_Q)$ of complex codimension one and  for each choice of the appropriate multipoles $\{w_k^{f_0} \}$,  the lifting of the deformation $f_t $  to the multipole spaces  is unique. $\Box$
 \end{itemize} 
 \end{theorem}

 Similarly, we get

\begin{theorem} Let $Q(x, y, z)$ be an irreducible real  quadratic form. Let
$\mathcal S_Q^\mathbb R: =  \{Q = 1 \}$ be a real quadratic surface,  equipped a  measure $dm$ for which any polynomial in $x, y$, and $z$ is an $L_2$-integrable function on the surface. Let  $f$ be an $L_2$-integrable continuous function on $\mathcal S_Q^\mathbb R$.  
 
Then there exists sequence  of multipoles 
 $\{w_k^f \in \overline {\mathcal M}(k)\}_{0 \leq k < +\infty}$   such that:
 \begin{itemize}
 \item the sequence  gives rise, via the maps $\{\overline \Phi_Q^\mathbb R\}$,  to real polynomials $P^f_d(x,y,z) = \sum_{k = 0}^d \overline \Phi_Q^\mathbb R(w_k^f)$, where the mutually orthogonal polynomials $\{\overline \Phi_Q^\mathbb R(w_k^f)\}$ are uniquely determined by $f$;
 \item As $d \rightarrow \infty$,   the polynomials $\{P_d^f\}$ converge to $f$ in the space 
 $L_2(\mathcal S_Q^\mathbb R)$;
\item The multipoles  $\{w_k^f \}$ satisfy  property $(57)$ from Corollary $20$.
\end{itemize}
When $Q$ is positive-definite,  
\begin{itemize}
\item the  mutually orthogonal polynomials $\{\overline \Phi_Q^\mathbb R(w_k^f) \in Har_Q(k; \mathbb R)\}$;
\item  each multipole $w_k^f $ is uniquely determined by $f$ ;
\item  for any  continuous deformation 
$\{f_t \in C(\mathcal S_Q^\mathbb R)\}_{ 0 \leq t \leq 1}$,
 of the function $f = f_0$, there  exists a unique continuous deformation  $\{w_k^{f_t} \in \overline {\mathcal M}(k)\}$ of the $f_t$-representing  multipoles, such that $\{w_k^{f_0} = w_k^f\}$.   $\Box$
 \end{itemize}
 \end{theorem}

\bigskip

{\it Acknowledgments.} I am grateful to Jeff Weeks for  introducing me to the subject. My conversations with Michael Shubin about the ``linear aspects" of this investigation were equally enlightening. I also would like to thank Blaine Lawson for explaining to me a few facts and constructions that turned out to be very helpful.

\section{References}

[B] Bennett, C. et al., {\it First year  Wilkinson Microwave Anisotropy Probe (WMAP 1) observations: preliminary maps and basic results}, Astrophysical Journal Suppliment Series 148 (2003), 1-27.

[Ch] Chow,  W.-L., {\it On the equivalence classes of cycles in an algebraic variety}, Ann. of Math. 64 (1956), 450-479.

[CHS] Copi C.J., Huterer D., Starkman, G.D., {\it Multipole vectors---a new representation of the CMB sky and evidence for statistical anisotropy or non-Gaussianity at} $2 \leq l \leq 8$, to appear in Phys. Rev. D. (astro-ph/0310511).

[CH] Courant,  R.,  Hilbert, D., {\it Methods of Mathematical Physics}, v.1, Interscience Publishers, 1953,  514-521.

[D] Dennis M. R., {\it Canonical representation of spherical functions: Sylvester's theorem, Maxwell's multipoles and Majorana's sphere} (arXiv:math-ph/0408046 v1).

[DT] Dold, A., Thom, R., {\it Quasifaserungen und unendliche symmetrische producte}, Ann. of Math. (2) 67 (1956), 230-281. 

[EHGL] Erisen, H.K., Banday, A.J.,  G\'{o}rski, K.M.,  and Lilje, P.B., {\it Asymmetries in the cosmic microwave background anisotropy field}, Astrophysics. J. 605 (2004) 14-20, (arXiv:astro-ph/0307507).

[G] Gamelin, T. W., {\it Uniform Algebras}, Prentice-Hall, Englewood Cliffs, N.J., 1969.

[GuR] Gunning, R.C., Rossi, H., {\it Analytic Functions of Several Complex Variables}, Prentice-Hall, 1965.

[H] Hartshorne, R., {\it Algebraic Geometry}, Springer-Verlag , 1983.

[Ha] Hatcher, A., {\it Algebraic Topology}, Cambridge University Press, 2002.

[K] Katz., G., {\it How Tangents Solve Algebraic Equations, or a Remarkable Geometry of Discriminant Varieties}, Expositiones Mathematicae 21 (2003) 219-261.

[KW] Katz, G. Weeks, J., {\it Polynomial Interpretation of Multipole Vectors}, Phys. Rev. D.
70, 063527 (2004)  (arXiv:astro-ph/0405631).

[L] Lachi\`{e}ze-Rey, M., {\it Harmonic projection and multipole vectors}, preprint (arXiv:astro-ph/0409081). 

[M] Maxwell, J.C., 1891 {\it A Treatise on Electricity and Magnetism}, v.1, (3rd edition) Clarendon Press, Oxford, reprinted by Dover, 1954.

[N] Narasimhan, R., {\it Introduction to the Theory of Analytic Spaces}, Lecture Notes in Mathematics 25 (1966), Springer-Verlag.

[Sh] Shubin, M. A., {\it Pseudo-differential Operators and Spectral Theory}, Nauka, Moscow, 1978.

[S] Sylvester, J.J.,  {\it Note on Spherical Harmonics}, Philosophical  Magazine, v. 2m, 1876, 291-307 \& 400.
  
[S1] Sylvester, J.J., 400. {\it Collected Mathematical Papers}, v.3, 37-51, Cambridge University Press, Cambridge, 1909.

[TOH]  Tegmark, M., de Oliveira-Costa and Hamilton A.J.S., {\it A high resolution foreground cleaned CMB map from WMAP}, Phys. Rev. D. 68 (2003) 123523, (arXiv:astro-ph/0302496).

[W] Weeks, J., {\it Maxwell's Multipole Vectors and the CMB}, preprint (arXiv:astro-ph/ 0412231).

\bigskip

Department of Mathematics, Brandeis University, Waltham, MA 02454\hfill\break
e-mail: {\it gabrielkatz@rcn.com}

 \end{document}